# Modeling Heavy-Gas Dispersion in Air with Two-Layer Shallow Water Equations

**Alexandre Chiapolino [1,*], Sébastien Courtiaud [2], Emmanuel Lapébie [2] and Richard Saurel [1,3]**

[1] Scientific Research and Numerical Simulation (RS2N), 371 chemin de Gaumin, 83640 Saint-Zacharie, France ; alexandre.chiapolino@rs2n.eu

[2] Commissariat à l'énergie atomique et aux énergies alternatives, Direction des applications militaires (CEA-DAM), CEA-Gramat, F-46500 Gramat, France ; sebastien.courtiaud@cea.fr ; emmanuel.lapebie@cea.fr

[3] Aix-Marseille Univ, French National Centre for Scientific Research (CNRS), Centrale Marseille, Laboratory of Mechanics and Acoustics (LMA), 4 impasse Nikola Tesla, 13013 Marseille, France ; richard.saurel@univ-amu.fr

* Correspondence: alexandre.chiapolino@rs2n.eu

**Abstract:** Computation of gas dispersal in urban places or hilly grounds requires a large amount of computer time when addressed with conventional multidimensional models. Those are usually based on two-phase flow or Navier-Stokes equations. Different classes of simplified models exist. Among them, two-layer shallow water models are interesting to address large-scale dispersion. Indeed, compared to conventional multidimensional approaches, 2D simulations are carried out to mimic 3D effects. The computational gain in CPU time is consequently expected to be tremendous. However, such models involve at least three fundamental difficulties. The first one is related to the lack of hyperbolicity of most existing formulations, yielding serious consequences regarding wave propagation. The second is related to the non-conservative terms in the momentum equations. Those terms account for interactions between fluid layers. Recently, these two difficulties have been addressed in Chiapolino and Saurel (2018) and an unconditional hyperbolic model has been proposed along with a Harten-Lax-van Leer (HLL) type Riemann solver dealing with the non-conservative terms. In the same reference, numerical experiments showed robustness and accuracy of the formulation. In the present paper, a third difficulty is addressed. It consists of the determination of appropriate drag effect formulation. Such effects also account for interactions between fluid layers and become of particular importance when dealing with heavy-gas dispersion. With this aim, the model is compared to laboratory experiments in the context of heavy gas dispersal in quiescent air. It is shown that the model accurately reproduces experimental results thanks to an appropriate drag force correlation. This function expresses drag effects between the heavy and light gas layers. It is determined thanks to various experimental configurations of dam-break test problems.

---

## 1. Introduction

### *1.1. Heavy-Gas Dispersal*

Gas dispersal may occur in many urban places, industrial plants and natural environments. Knowledge of the flow field is important in such contexts because dispersal may carry pollutants or hazardous chemical species. For instance, certain industrial processes involve storing hazardous liquefied gases, which may be released accidentally. The release of these heavier-than-air gases is of serious interest regarding risk assessment.

From a numerical point of view, one of the difficulties is to address long-time and large-scale computations, for instance, at the scale of a city, while providing accurate results with reasonable CPU time. Multiple classes of dense-gas dispersion models are available in the literature, such as integral or box models [1,2], intermediate (SLAB-type) models [3], shallow water models [4], and multidimensional computational fluid dynamics (CFD) models [5]. Relevant literature on this subject can be found in the PhD thesis of Hankin, 1997 [4].

Integral or box models typically regard an instantaneously released cloud of dense gas as a uniform cylinder. Those are the simplest form of heavy-gas dispersion models. However, they do not have flexibility to characterize a cloud influenced by a complex terrain. They are limited by the assumption that the cloud remains as a uniform cylinder or ellipse. Such an assumption is obviously inappropriate if the cloud is channeled by topographical features into low-lying areas for example.

Multidimensional computations based on two-phase flow or Navier–Stokes equations are obviously much more complex and are capable of accounting for a complex terrain. However, their numerical resolution demands a large amount of computer time. Illustrative 3D results are provided in Hank et al., 2012, 2014 [5,6], for example.

Between the simplicity of integral or box models and the complexity of multidimensional CFD models, intermediate approaches have been proposed. Examples are SLAB-type models, where a Lagrangian description is used to track gas bubbles of a plume. Another approach is based on shallow water models [4].

Results obtained with the commonly used SLAB-type models are provided in Hanna et al., 2008 [7], for instance. Although such models admit more complex cloud geometries than integral models, they also suffer from similar drawbacks as they are inappropriate for describing a cloud influenced by complex topographies [4].

The shallow water approach has a number of advantages over other simplified models as thoroughly argued in Hankin, 1997 [4]. The main advantage of such approach is that it is able to account for complex topographies such as buildings, valleys and mountain ranges.

As only the height and speed of the dense gas layer are of interest, a single-layer shallow water model is attractive. Nevertheless, in the present heavy-gas dispersion context, the single-layer shallow water approach is inaccurate as shown in the numerical experiments of Chiapolino and Saurel, 2018 [8]. The lack of accuracy is due to the low-density ratio between the two gases. Interactions between the ambient gas and the dense gas layer are indeed important and are to be accounted for [4,8].

In this direction, the present paper presents a novel approach. The dispersal of the heavy gas in air is treated with the recent hyperbolic two-layer shallow water model of Chiapolino and Saurel, 2018 [8]. Interactions between the two gases are considered through the non-conservative terms of the momentum equations and through drag effects, occurring at the interface between the dense gas and the ambient air.

Turbulence effects are considered at the interface through the corresponding interfacial area. However, determination of proper interfacial aera is a hard task. In the present contribution, it is adjusted to address heavy gas dispersion into a quiescent atmosphere and on a flat ground. This step is indeed necessary to validate the hyperbolic two-layer shallow water model in its simple form, i.e., considering only wave propagation and drag effects.

More complex and realistic situations will be part of future investigations, such as the introduction of multidimensional and topography effects as well as turbulence. Relevant literature on the last subject is, for instance, the works of Teshukov, 2007 [9] and Richard and Gavrilyuk, 2012 [10] where a model for turbulent shear shallow water flows is derived. In this latter contribution, the rapidly varied flows studied by the authors are characterized by the presence of a turbulent structure called roller in which the turbulent energy dissipation plays a major part. This work has been continued by the same authors in [11–14], to cite a few.

However, the aforementioned references address turbulent structures with a one-layer shallow water type model. Adaptation to two-layer equations is out of the scope of the present paper but shall also be considered in future investigations.

The present approach based on two-layer shallow water equations is very different from usual atmospheric models that are unsuitable for describing a cloud of dense gas influenced by urban structures or hilly grounds. Topography effects are omitted in the present contribution for the sake of simplicity but shall be addressed in future works as well. Indeed, even in the absence of topography effects or agitated atmospheric conditions, two-layer models need modelling and validation efforts, as detailed a bit further.

It is therefore important to address first heavy-gas dispersion on a flat ground and in a quiescent atmosphere. Such task is addressed in the present paper. It is shown that the two-layer approach has nice predicting abilities. A background of two-layer shallow water models is presented hereafter.

### *1.2. Two-Layer Shallow Water Approach*

The two-layer shallow water strategy is indeed an interesting candidate as it allows addressing 2D simulations to mimic 3D flows. The computational gain in CPU time is thereby expected to be tremendous compared to conventional multi-fluid approaches, see for example Saurel and Pantano, 2018 [15] for a review of diffuse interface capturing methods. In Chiapolino and Saurel, 2018 [8], CPU time saving factor of the order of million is reported thanks to a two-layer shallow water approach.

Two-layer (and multi-layer) shallow water models are particularly useful in some limit cases of multi-fluid and variable density flows separated by nearly horizontal interfaces. These models govern the dynamics of incompressible fluids spreading under gravity effects. It can be, for example:

- Flows of the same liquid but at different temperatures, resulting in density differences, such situations being typical of oceanic flows.
- Flows of two liquids of different densities.
- Flows of two gases evolving at low Mach number.

In this contribution, this last situation is of particular interest. As mentioned earlier, heavy-gas dispersal in urban places or hilly grounds motivates the present work.

The two-layer approach is also helpful when the height of one of the phases is arbitrarily small as there is no need to spatially resolve it. Thereby, no numerical diffusion of the nearly horizontal interface is present, and no interface tracking is needed. Such models can deal with topography effects or obstacles (not treated in the present paper). However, there are obviously some limitations with this approach:

- The vertical velocity component is neglected.
- The velocity is assumed uniform in cross sections of each layer.

Such types of modeling also involve serious difficulties. Indeed, most models are not hyperbolic, this issue having serious consequences both for wave propagation, which becomes ill-posed, and for the design of numerical methods. A second serious difficulty appears as non-conservative terms are present in the momentum equations. Recently, those difficulties have been addressed in Chiapolino and Saurel, 2018 [8] and solutions have been proposed. A third difficulty appears in the modeling of drag effects between fluid layers. This issue is addressed in the present contribution with the help of dedicated laboratory dam-break experiments.

In the frame of averaged (or homogenized) equations in fluid mechanics, the issue related to the lack of hyperbolicity appears in different types of models, such as those of non-equilibrium two-phase flows. Only a few models seem well-posed with this respect (Marble, 1963 [16], Baer and Nunziato, 1986 [17], Saurel et al., 2017 [18]). There are mainly two types of remedy to cure this issue:

- Consider compressibility of the phases and deal with pressure relaxation (Lallemand and Saurel, 2000 [19]). This approach involves sound propagation in the phases and is particularly efficient in many situations. It has been adopted in [17,18].
- Consider turbulent effects in the phases, as they result in the appearance of a "turbulent sound speed" (Forestier et al., 1997 [20], Saurel et al., 2003 [21], Lhuillier et al., 2013 [22]). In the frame of shallow water flows, these effects have been studied in Richard and Gavrilyuk, 2012 [10] and Gavrilyuk et al., 2016 [14].

In the present work, the first method is adopted, and the fluids are considered weakly compressible. The resulting model is unconditionally hyperbolic and in the limit of stiff pressure relaxation, the conventional (non-hyperbolic) two-layer model is recovered. This approach is reminiscent of the model of Abgrall and Karni, 2009 [23] except that extra pressure terms are present in the momentum equations of the new formulation. It also gives another interpretation of the relaxation approach, now based on compressibility and pressure effects.

The second issue, related to the presence of non-conservative terms in the momentum equations, has also been addressed in [8]. By examining the Riemann problem structure, it appears that local constants are present, at locations where the derivative of the Heaviside function emerges. Consequently, the non-conservative products become well-defined and a locally conservative form is obtained and used in the frame of a Harten-Lax-van Leer (HLL) type Riemann solver.

The accuracy and robustness of the new HLL-type solver has been checked against results of Abgrall and Karni, 2009 [23] as well as results obtained with a flow solver based on the VFRoe method of Gallouet and Masella, 1996 [24] as it is able to deal, to some extent, with both conservative and non-conservative systems.

The theoretical formulation of Chiapolino and Saurel, 2018 [8] and its associated HLL-type solver make consequently an interesting strategy to address large-scale gas dispersal. The 1D two-layer model has been confronted in [8] to 2D diffuse interface computations [15,25–27] and showed good agreement. However, a parameter related to the drag force was needed to match 1D and 2D averaged results in the vertical direction. In the present paper, the model is confronted to laboratory experiments, specifically designed for this aim. A drag function is proposed, with three parameters determined explicitly, enabling accurate matching of computed 1D results and experiments, in a broad range of initial conditions of two-layer dam-break test problems.

This paper is organized as follows. The two-layer model with its hyperbolic formulation is presented in Section 2. Pressure and velocity relaxation (drag effects) are introduced as well. The numerical treatment is summarized in Section 3. The experimental facility used for dam-break-type tests is presented in Section 4. Comparisons between 1D computed results and laboratory experiments are presented in Section 5, enabling building of an appropriate drag function between heavy and light fluid layers. It is shown that the hyperbolic model accurately reproduces experimental results at the price of few parameters related to the drag function. Conclusions are given in Section 6.

**2. Hyperbolic Two-Layer Shallow Water Model**

The model of Chiapolino and Saurel, 2018 [8] reads:

$$\begin{cases}\dfrac{\partial h_1}{\partial t}+u_1\dfrac{\partial h_1}{\partial x}=\mu\dfrac{\left(p_1-p_0-\rho_2 g h_2\right)}{\rho_1 c_1^2}, & \mu\to+\infty,\\ \dfrac{\partial\left(h_1\rho_1\right)}{\partial t}+\dfrac{\partial\left(h_1\rho_1 u_1\right)}{\partial x}=0, & \\ \dfrac{\partial\left(h_1\rho_1 u_1\right)}{\partial t}+\dfrac{\partial\left(h_1\rho_1 u_1^2+h_1 p_1\left(\rho_1,\rho_2,h_2\right)+\frac{1}{2}\rho_1 g h_1^2\right)}{\partial x}=\rho_2 g h_2\dfrac{\partial h_1}{\partial x}+p_0\dfrac{\partial h_1}{\partial x}, & \\ \dfrac{\partial h_2}{\partial t}+u_2\dfrac{\partial h_2}{\partial x}=\mu\dfrac{\left(p_2-p_0\right)}{\rho_2 c_2^2}, & \\ \dfrac{\partial\left(h_2\rho_2\right)}{\partial t}+\dfrac{\partial\left(h_2\rho_2 u_2\right)}{\partial x}=0, & \\ \dfrac{\partial\left(h_2\rho_2 u_2\right)}{\partial t}+\dfrac{\partial\left(h_2\rho_2 u_2^2+h_2 p_2\left(\rho_2\right)+\frac{1}{2}\rho_2 g h_2^2\right)}{\partial x}=-\rho_2 g h_2\dfrac{\partial h_1}{\partial x}+p_0\dfrac{\partial h_2}{\partial x}. & \end{cases} \quad (1)$$

In the present paper, index 1 denotes the heaviest fluid (lower layer) and index 2 the lightest one (upper layer). The notations are conventional, $h_1$ and $h_2$ denote the heights of the two layers, $\rho_1$ and $\rho_2$ represent the densities of the fluids, $u_1$ and $u_2$ denote the fluid velocities, averaged in each layer and $g$ represents the gravity constant. $p_1$ and $p_2$ denote the thermodynamic pressure of the fluids, given by barotropic (and convex) equations of state. An example of such an equation of state (EOS) is:

$$\begin{cases}p_1=p_0+\rho_2 g h_2+c_1^2\left(\rho_1-\rho_1^{(0)}\right),\\ p_2=p_0+c_2^2\left(\rho_2-\rho_2^{(0)}\right).\end{cases} \quad (2)$$

Other options, such as Tait EOS are possible. As shown in [8], the choice of the EOS is not important, only the related fluid sound speeds $c_1$ and $c_2$ have importance. $p_0$ denotes the (constant) atmospheric pressure. Hydrostatic effects are then considered through the EOS of the first fluid, i.e., the heaviest layer.

Topography effects have been omitted for the sake of simplicity as well as friction with the bottom and between layers. As shown in [8], System (1) recovers the conventional but conditionally hyperbolic model of Ovsyannikov, 1979 [28] in the stiff pressure relaxation limit $\left(\mu\to+\infty\right)$. Such conventional model has been examined in Abgrall and Karni, 2009 [23], Kurganov and Petrova, 2009 [29] and Monjarret, 2015 [30] and appeared hyperbolic for low velocity drift only:

$$\left(u_1-u_2\right)^2<\left(h_1-h_2\right)g\left(1-\frac{\rho_2}{\rho_1}\right).$$

In the present approach, pressure non-equilibrium effects result in an unconditionally hyperbolic formulation. Such hyperbolic system has been developed in Abgrall and Karni, 2009 [23] and an extended version in Chiapolino and Saurel, 2018 [8].

Compared to the conventional two-layer shallow water model (Ovsyannikov, 1979 [28]), two equations have been added and express the transport of the heights of the fluid layers that are assumed to vary as a function of pressure differentials.

The pressure relaxation parameter $\mu$ is related to the sound speeds and heights of fluid layers. It controls the rate at which pressure equilibrium is reached. In most situations, the relaxation time $\tau$ is of the order of 1/100 s meaning that the relaxation parameter $\mu$ is large: $\mu \simeq \max\left(\tau_1^{-1}, \tau_2^{-1}\right)$ or alternatively $\mu \simeq \min\left(c_1, c_2\right)$. Details are provided in Chiapolino and Saurel, 2018 [8].

In practical computations, the relaxation time $\tau$ will be assumed of the same order as the computational time step and stiff pressure relaxation will be done at the end of each time step. Therefore, there is no need for a precise knowledge of the pressure relaxation parameter $\mu$.

System (1) is hyperbolic with wave speeds $u_k$, $u_k + \sqrt{c_k^2 + \frac{1}{2} g h_k}$, $u_k - \sqrt{c_k^2 + \frac{1}{2} g h_k}$ with $k = \{1, 2\}$. This model is however relevant with respect to the physics expressed in the conventional but non-hyperbolic system of Ovsyannikov, 1979 [28], as it tends to the same equations when pressure relaxation is stiff, see [8] for details.

Consequently, as System (1) is hyperbolic, it is a good candidate to numerically approximate the conventional model with a two-step procedure:

- Resolution of the hyperbolic step, i.e., the resolution of System (1) in the absence of source terms. An HLL-type Riemann solver has been developed in Chiapolino and Saurel, 2018 [8] and is summarized in the following.
- Stiff pressure relaxation towards the hydrostatic and atmospheric pressures and reset of the heights. This step makes the solution tend to the one of the conventional models [28].

In the present context, pressure relaxation is particularly simple. It just consists of resetting the heights of the layers $h_k \to h_k^*$ and is independent of the equation of state [8]. For the lightest fluid, thanks to mass conservation, the relaxed state reads:

$$h_2^* = \frac{h_2 \rho_2}{\rho_2^{(0)}}. \tag{3}$$

For the heaviest fluid, mass conservation and hydrostatic effects result in the relaxed state:

$$h_1^* = \frac{h_1 \rho_1}{\rho_1^{(0)} + \frac{\rho_2^* g h_2^*}{c_1^2}}. \tag{4}$$

In these relations, $\rho_k^{(0)}$ denotes the constant density of fluid $k$ at atmospheric pressure and superscript * denotes the relaxed solution state.

The hyperbolic two-layer shallow water model is able to recover the solution of the conventional but non-hyperbolic two-layer system. Besides, the single-layer Saint-Venant solution is recovered as well when the flow conditions are appropriate. Such situation happens either when the height of the heaviest layer is insignificant $h_1 \to 0$ or when the effects of the surrounding fluid are negligible. This behavior appears when the density ratio between the lightest fluid and the heaviest one is very small: $r = \frac{\rho_2}{\rho_1} \ll 1$. Examples are provided in Chiapolino and Saurel, 2018 [8]. In the same reference, the importance of the two-layer model when the density ratio $r = \rho_2 / \rho_1$ is arbitrary, is illustrated. Indeed, the two-layer system can deal with interactions between fluids unlike the conventional one-fluid Saint-Venant model.

Consequently, the one-layer model yields inaccurate results when the order of magnitude of the density ratio $r = \rho_2 / \rho_1$ is about $[0.2 - 0.5]$, see [8] for illustrative results and further discussions. In the present paper, such situation is of particular interest as the flow involves only gases.

The two-layer model is then necessary. The present model has been confronted in [8] to 2D diffuse interface computations [15,25–27] and showed good agreement after the introduction of drag effects. Pressure (or "acoustic") drag is considered only and is modeled through velocity relaxation terms that appear in the right-hand side of the momentum equations. For non-viscous fluids, the "acoustic" drag force can be written as:

$$F = \frac{Z_1 Z_2}{Z_1 + Z_2} A_I \left(u_2 - u_1\right). \tag{5}$$

$Z_k = \rho_k c_k$ denotes the acoustic impedance of fluid $k$ and $A_I$ denotes the specific interfacial area. This acoustic drag effect modeling was developed in Saurel et al., 2003 [21] and Chinnayya et al., 2004 [31]. It is obtained by local interfacial pressure integration over a piece of interface. The interfacial pressure is estimated through local approximate Riemann solver for the Euler equations of gas dynamics. Let us mention that upon integration over height, the specific interfacial area $A_I$ becomes dimensionless. In Chiapolino and Saurel, 2018 [8], $A_I$ was considered as a constant parameter. In Section 5, a correlation based on experimental results will be built. Beforehand, let us present the HLL-type Riemann solver used during the hyperbolic step.

## 3. HLL-Type Riemann Solver and Godunov-Type Method

A Riemann solver based on the Rankine-Hugoniot relations, such as the HLL solver of Harten et al., 1983 [32], has been developed in [8]. This HLL-type solver is simple, accurate and robust. It is summarized in the following.

For each fluid $k$, the two extreme, left- and right-facing waves ($S_{L,k}$ and $S_{R,k}$ respectively) are approximated following Davis, 1988 [33] as:

$$\begin{cases} S_{L,k} = \min\left( u_{L,k} - \sqrt{c_{L,k}^2 + \frac{1}{2} g h_{L,k}} \quad , \quad u_{R,k} - \sqrt{c_{R,k}^2 + \frac{1}{2} g h_{R,k}} \right), \\ S_{R,k} = \max\left( u_{L,k} + \sqrt{c_{L,k}^2 + \frac{1}{2} g h_{L,k}} \quad , \quad u_{R,k} + \sqrt{c_{R,k}^2 + \frac{1}{2} g h_{R,k}} \right), \end{cases} \tag{6}$$

with $k = \{1, 2\}$. The indexes L and R denote respectively the left and right states at a given cell boundary. The two most extreme waves ($S_L$ and $S_R$), used in the Riemann problem solution, are then approximated as:

$$S_L = \min\left(S_{L,1}, S_{L,2}\right), \ S_R = \max\left(S_{R,1}, S_{R,2}\right). \tag{7}$$

The two contact waves $u_1$ and $u_2$ are considered as well for the transport of the heights $h_1$ and $h_2$. Their exact Riemann problem solutions are straightforward:

$$\begin{cases} h_1^*\left(\frac{x}{t} < u_1^*\right) = h_{1,L}, \ h_1^*\left(\frac{x}{t} > u_1^*\right) = h_{1,R}, \\ h_2^*\left(\frac{x}{t} < u_2^*\right) = h_{2,L}, \ h_2^*\left(\frac{x}{t} > u_2^*\right) = h_{2,R}. \end{cases} \tag{8}$$

In the present section, the superscript $*$ denotes the Riemann problem solution. Relations (8) indicate that the non-conservative terms have contributions between the two extreme waves $S_L$ and $S_R$, at points where $h_1$ and $h_2$ are discontinuous. More precisely, only the discontinuity in $\rho_2 g h_2 \frac{\partial h_1}{\partial x}$ needs attention, as the non-conservative terms involving the atmospheric pressure (considered constant) transform to fluxes: $p_0 \frac{\partial h_k}{\partial x} = \frac{\partial (p_0 h_k)}{\partial x}$.

An analysis of the Riemann problem structure is provided in Chiapolino and Saurel, 2018 [8] and reveals that local constants appear precisely at locations where the non-conservative products need definition. Thanks to these local constants, a locally conservative formulation of the equations is obtained. Consequently, the non-conservative terms do not cause numerical difficulties during the resolution of the Riemann problem, this asset being significant.

Besides, to simplify the algorithm, a single solution state is considered for the apparent densities $(\rho_k h_k)^*$ instead of the two $(\rho_k h_k)_L^*$ and $(\rho_k h_k)_R^*$ in the same spirit as in the HLL solver for the Euler equations:

$$U_{k,mass}^* = (\rho_k h_k)^* = \frac{(\rho_k h_k)_R (u_{k,R} - S_R) - (\rho_k h_k)_L (u_{k,L} - S_L)}{S_L - S_R}. \quad (9)$$

Such approximation relying on two waves instead of three remains accurate in the present context, as the intermediate wave vanishes at the end of the pressure relaxation step. It is present only during the hyperbolic step, not in the asymptotic system resulting of pressure relaxation. Thanks to this approximation, the momentum equations become locally:

$$\begin{cases} \frac{\partial (h_1 \rho_1 u_1)}{\partial t} + \frac{\partial \left( h_1 \rho_1 u_1^2 + h_1 \left( p_1(\rho_1, \rho_2, h_2) - p_0 \right) + \frac{1}{2} \rho_1 g h_1^2 - g (\rho_2 h_2)^* h_1 \right)}{\partial x} = 0 \quad \Leftrightarrow \quad \frac{\partial U_{1,mom}}{\partial t} + \frac{\partial F_{1,mom}}{\partial x} = 0, \\ \frac{\partial (h_2 \rho_2 u_2)}{\partial t} + \frac{\partial \left( h_2 \rho_2 u_2^2 + h_2 \left( p_2(\rho_2) - p_0 \right) + \frac{1}{2} \rho_2 g h_2^2 + g (\rho_2 h_2)^* h_1 \right)}{\partial x} = 0 \quad \Leftrightarrow \quad \frac{\partial U_{2,mom}}{\partial t} + \frac{\partial F_{2,mom}}{\partial x} = 0. \end{cases} \quad (10)$$

Solutions in terms of conservative variables $U_k^*$ and fluxes $F_k^*$ are then computed by the HLL approximation:

$$U_k^* = \frac{F_{k,R} - F_{k,L} - S_R U_{k,R} + S_L U_{k,L}}{S_L - S_R}, \quad F_k^* = \frac{F_{k,R} S_L - F_{k,L} S_R + S_L S_R (U_{k,L} - U_{k,R})}{S_L - S_R}. \quad (11)$$

The speeds of the fluids are given by:

$$u_k^* = \frac{U_{k,mom}^*}{U_{k,mass}^*} = \frac{(h_k \rho_k u_k)^*}{(h_k \rho_k)^*}. \quad (12)$$

The various equations of System (1) are updated with a Godunov-type method (stable under the conventional CFL condition) as:

$$\begin{cases} h_{k,i}^{n+1} = h_{k,i}^n - \frac{\Delta t}{\Delta x} \left( (hu)_{k,i+1/2}^* - (hu)_{k,i-1/2}^* \right) + \frac{\Delta t}{\Delta x} h_{k,i}^n \left( u_{k,i+1/2}^* - u_{k,i-1/2}^* \right), \\ (h_k \rho_k)_i^{n+1} = (h_k \rho_k)_i^n - \frac{\Delta t}{\Delta x} \left( F_{k,mass,i+1/2}^* - F_{k,mass,i-1/2}^* \right), \\ (h_1 \rho_1 u_1)_i^{n+1} = (h_1 \rho_1 u_1)_i^n - \frac{\Delta t}{\Delta x} \left( F_{1,mom,i+1/2}^* - F_{1,mom,i-1/2}^* \right) + \frac{\Delta t}{\Delta x} h_{1,i}^n \left( -g \left[ (h_2 \rho_2)_{i+1/2}^* - (h_2 \rho_2)_{i-1/2}^* \right] \right), \\ (h_2 \rho_2 u_2)_i^{n+1} = (h_2 \rho_2 u_2)_i^n - \frac{\Delta t}{\Delta x} \left( F_{2,mom,i+1/2}^* - F_{2,mom,i-1/2}^* \right) + \frac{\Delta t}{\Delta x} h_{1,i}^n \left( g \left[ (h_2 \rho_2)_{i+1/2}^* - (h_2 \rho_2)_{i-1/2}^* \right] \right), \end{cases} \quad (13)$$

where n+1 and n denote two consecutive time steps and superscript $*$ denotes the Riemann problem solution. Indexes $i$ and $i \pm 1/2$ denote respectively the center of the current numerical cell and its corresponding boundaries. Note that the present integration of the non-conservative terms assumes $h_k$ constant and taken at cell center $h^n_{k,i}$.

The HLL-type solver and its associated Godunov-type scheme (13) provide accurate results, maintain height positivity even in extreme situations and are oscillation free as shown in [8]. In this reference, it is also shown that sound speeds $c_k$ influence computed results and particularly numerical diffusion.

For the sake of clarity, let us consider a limit situation where $h_1 \to 0$. System (1) then tends to the one-layer Saint-Venant equations, in the limit of stiff pressure relaxation. The wave speeds of the two-layer model (1) involve the effective sound speeds given by $\pm\sqrt{c_k^2 + \frac{1}{2}gh_k}$ while the single-layer wave speeds involve only $\pm\sqrt{gh}$. As shown in [8], when $c_k$ is significantly greater than $\sqrt{\frac{1}{2}gh_k}$ excessive numerical diffusion is present. However, computational results do converge towards the exact solution when the mesh is fine enough. It thus appears that large sound speeds are admissible but result in excessive numerical diffusion. Numerical experiments also reveal that the method becomes unstable when $c_k < \sqrt{\frac{1}{2}gh_k}$. They consequently suggest existence of a subcharacteristic condition. The following numerical sound speed has consequently been proposed in [8],

$$c_k = \sqrt{\theta_k \frac{1}{2} gh_k}, \quad \text{with } \theta_k > 1, \quad k = \{1, 2\}, \tag{14}$$

during the hyperbolic step only, i.e., during the resolution of the Riemann problem.

$\theta_k$ is a parameter that controls numerical diffusion. In order to unambiguously fulfill the above-mentioned subcharacteristic condition, $\theta_k$ must be greater than unit. According to numerical experiments, $\theta_k \in [2, 5]$ seems to be a fair choice as it is low enough to control numerical diffusion and large enough to ensure stability (see Chiapolino and Saurel, 2018 [8] for details).

## 4. Experimental Apparatus

Many experimental results dealing with fluid dispersal under gravity effects are provided in the literature and compared to shallow water computations. However, most of the experiments are carried out with one or several liquids, see for example Adduce et al., 2012 [34] and references therein. Consequently, a dedicated experimental setup was designed to produce and characterize an axisymmetric flow of dense gas, and thus obtain the necessary data for the validation of the model. The principle of the experiment is similar to the one used during the full-scale trials of Thorney Island [35] and consists of filling a cylindrical container with a dense gas and releasing it instantaneously (or very fast) by letting the container's walls fall down, thus creating a gravity-driven axisymmetric flow. The current tests were performed at small scale and pure krypton was used as the dense gas. In the present section, krypton will be specifically denoted with the index "Kr". The second gas is the environmental air and will be denoted specifically with the index "a".

### *4.1. Description of the Experimental Setup*

The experimental setup is depicted in Figure 1. It is made of a closed, 1.2-m large, cubic vessel. The vessel is not airtight and is only used to isolate the experiment from any external perturbation. Three of the lateral faces, the ceiling and the floor are transparent to have optical access. The latter is made of plexiglass whereas the other faces are made of glass. It is to be noted that the floor is supported by the aluminum frame of the vessel but also by two additional beams. The vessel is placed in an open room which means that the external conditions could not be controlled during the trials.

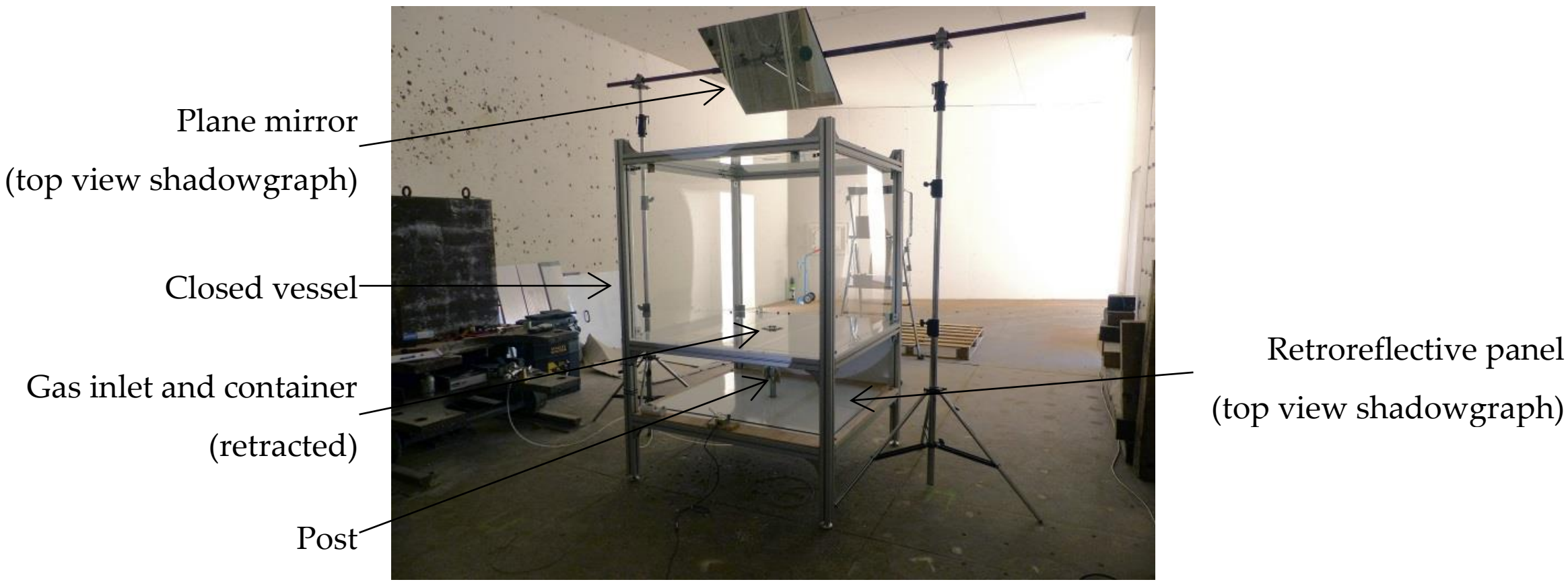

**Figure 1.** Global view of the experimental setup. The retro-reflective back panel was removed when the picture was taken.

The floor is equipped with a retractable metallic cylindrical container (shown in Figure 2) filled with krypton at the beginning of each trial. Its internal diameter is 99.6 mm and its height above the floor could be set to three different values: 6 cm, 11 cm and 16 cm, corresponding respectively to gas heights of 5 cm, 10 cm, and 15 cm. When retracted, the top of the container is flush with the floor. The container is installed on a post on which it could slide, and which includes the gas inlet. The container is initially maintained in position thanks to a pin, visible in Figure 2b, that is simply removed at the start of the experiment. In order to avoid the apparition of strong vibrations when the container reaches its lower position, an elastic damper is installed on the post. Besides, the rebound of the container is mitigated thanks to two sandbags placed at its bottom. The airtightness of the junction between the floor and the container is ensured by a rubber gasket and holes are made at the bottom of the container to avoid any suction effect during its descent. The gas inlet is 4 mm in diameter and a metallic part is put above it to break the krypton jet and ensure a homogeneous filling of the container from the bottom to the top with pure gas.

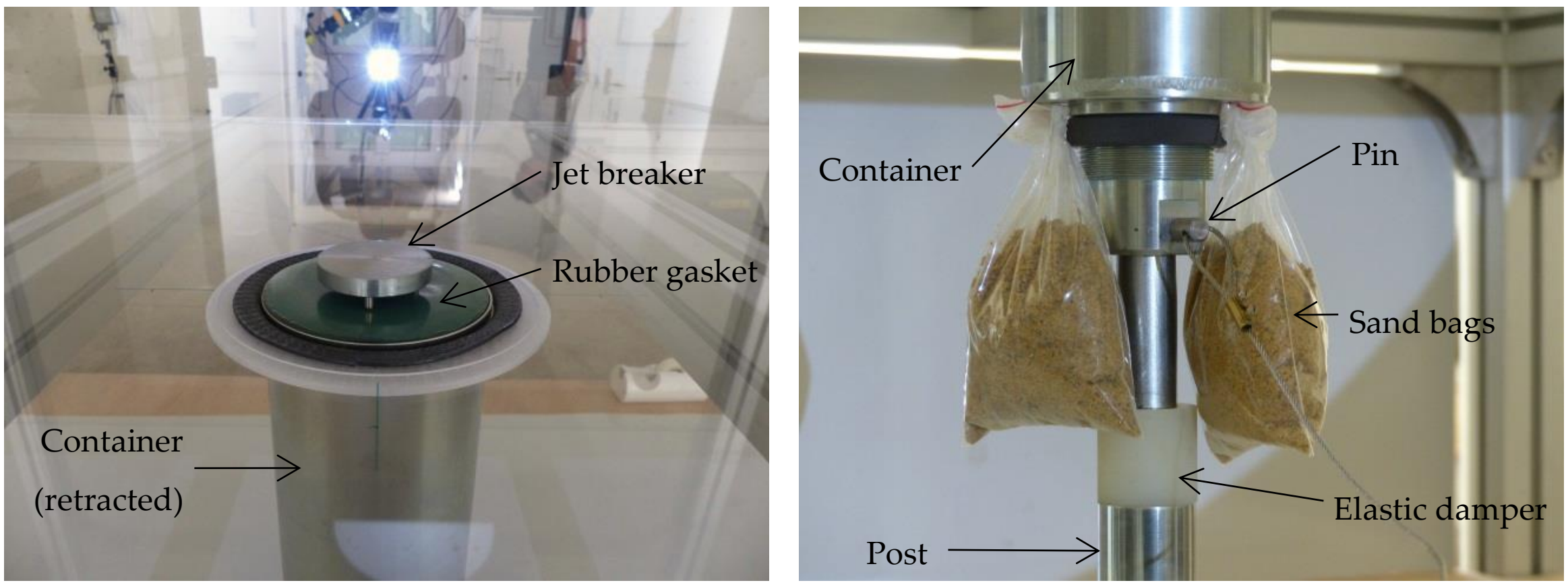

(**a**) View from above the floor (**b**) View from below the floor

**Figure 2.** Detailed view of the retractable, metallic, cylindrical container used to contain the dense gas at the initial state.

The resulting krypton cloud is visualized along two axes using Edgerton's, 1958, shadowgraphy technique [36]. A first camera looks at the cloud in a vertical transversal plane (front view) whereas a second one looks at the cloud from above (top view). To this end, the back panel is covered with a retro-reflective material, a similar panel is placed below the floor of the vessel and a plane mirror is positioned above the vessel. The lenses of the two cameras are equipped with small mirrors to direct the light on the same axis as the camera optical axis as shown in Figure 3. The cameras used are a Phantom V2511 for the front view and a Photron Fastcam SA-Z for the top view. The frame rate is 500 Hz for both cameras. The resolution of the shadowgraphs is 0.86 mm per pixel for the front view and 0.97 mm per pixel for the top view.

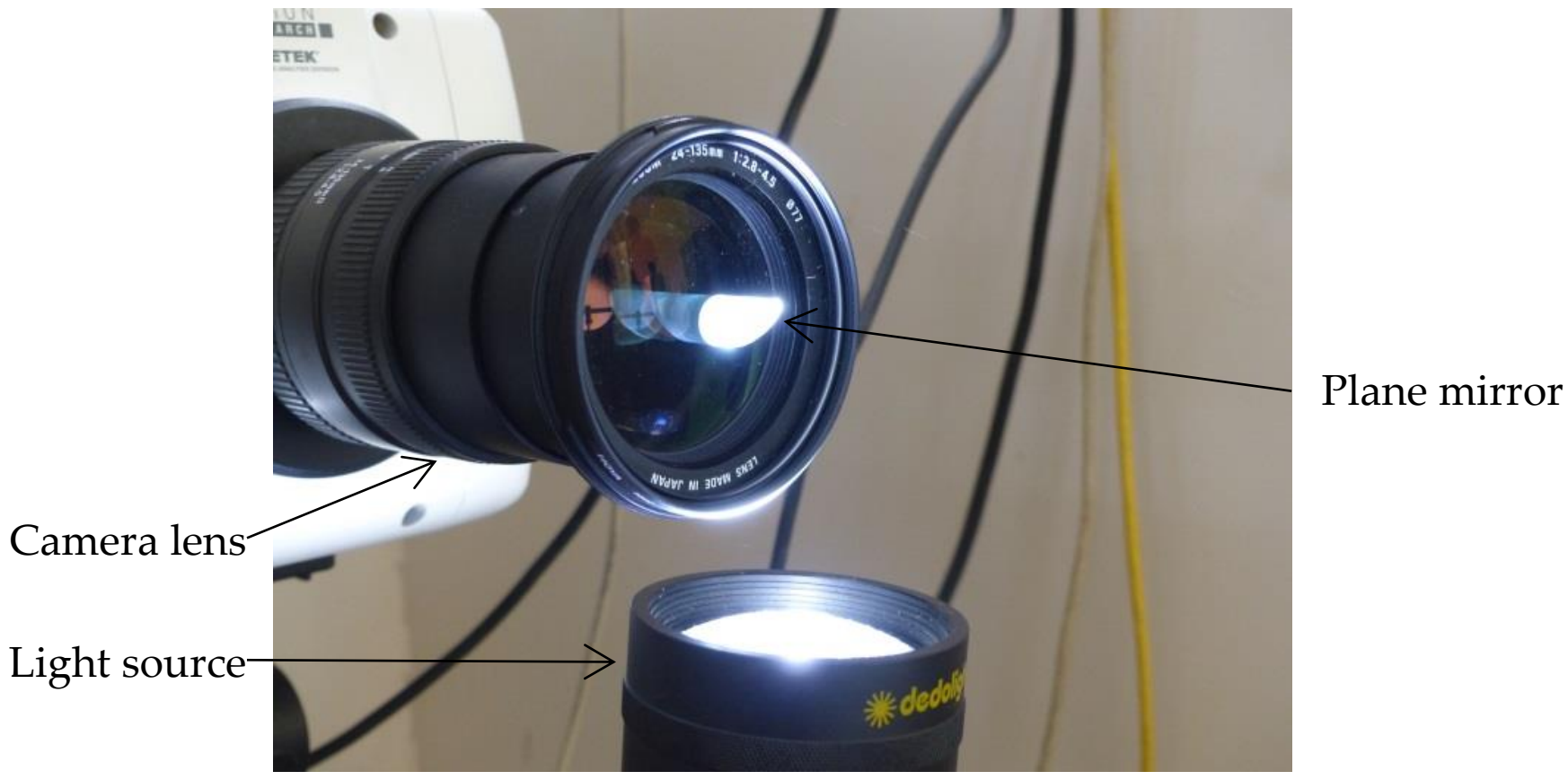

**Figure 3.** Detailed view of the lens of one of the cameras showing the small plane mirror used to direct the light on the same axis as the camera optical axis.

The amount of gas introduced inside the container is controlled with a ball flowmeter and a chronometer. Given the height $h_{Kr}$ of krypton to be introduced, and a target value of the volumetric flow rate $Q_{Kr}$ of krypton, the target value for the duration of the filling process $t_{Kr}$ is determined as,

$$t_{Kr} = \frac{\pi R_c^2 h_{Kr}}{Q_{Kr}}, \tag{15}$$

where $R_c$ is the internal radius of the container.

### 4.2. Course of Operations

Each test is performed as follows. First, the back panel of the vessel is removed, and the container is lowered. The bottle of krypton and the gas inlet are open, and the flowmeter is set to the target flow rate value. The gas inlet is closed, and the vessel is aerated. The container is then set to the desired height and the vessel is closed. After a few minutes (according to the estimate of Equation (15)), the container is filled with krypton. The experiment is started by removing the pin at the bottom of the container. Krypton is thereby released almost instantaneously and creates a cloud spreading on the floor. The time between the end of the filling process and the removal of the pin did not exceed a few seconds.

### 4.3. Processing of the Shadowgraphs

The analysis of the shadowgraphs and the extraction of the radius of the cloud have been performed automatically thanks to a dedicated tool written in Python and making use of the image processing library Scikit-Image [37]. Since the images of the top and front view differ, the processing method has been adapted to each view.

In the top view images, the cloud of krypton has a circular shape. Consequently, the processing method consists of determining a circle which fits the contours of the cloud. The processing steps are illustrated in Figure 4. First, the background of the image is subtracted (Figure 4b). The background image is calculated by taking the average of several dozens of pictures that are taken between the moment when the pin leaves the field of view and the moment when the container reaches its lowest position. The contrasts of the picture are then enhanced using a Contrast Limited Adaptive Histogram Equalization (CLAHE) technique (Figure 4c) which results in a clean picture of the cloud alone. Then two filters are applied to determine the contours of the cloud. The first one is a Frangi filter [38] (Figure 4d), which is well-adapted to circular patterns, and the second one is a Canny filter [39] (Figure 4e) which makes the image binary and facilitates the fitting step. The last step consists of determining the circle that best fits the contours of the cloud thanks to the Hough transform technique (Figure 4f). The extraction of the contours is limited by the initial contrast of the images and thus by the thickness of the cloud. As it expands, the contrast decreases and the determination of the circle becomes more difficult. Another restriction comes from the field of view of the camera which is limited by the size of the plane mirror.

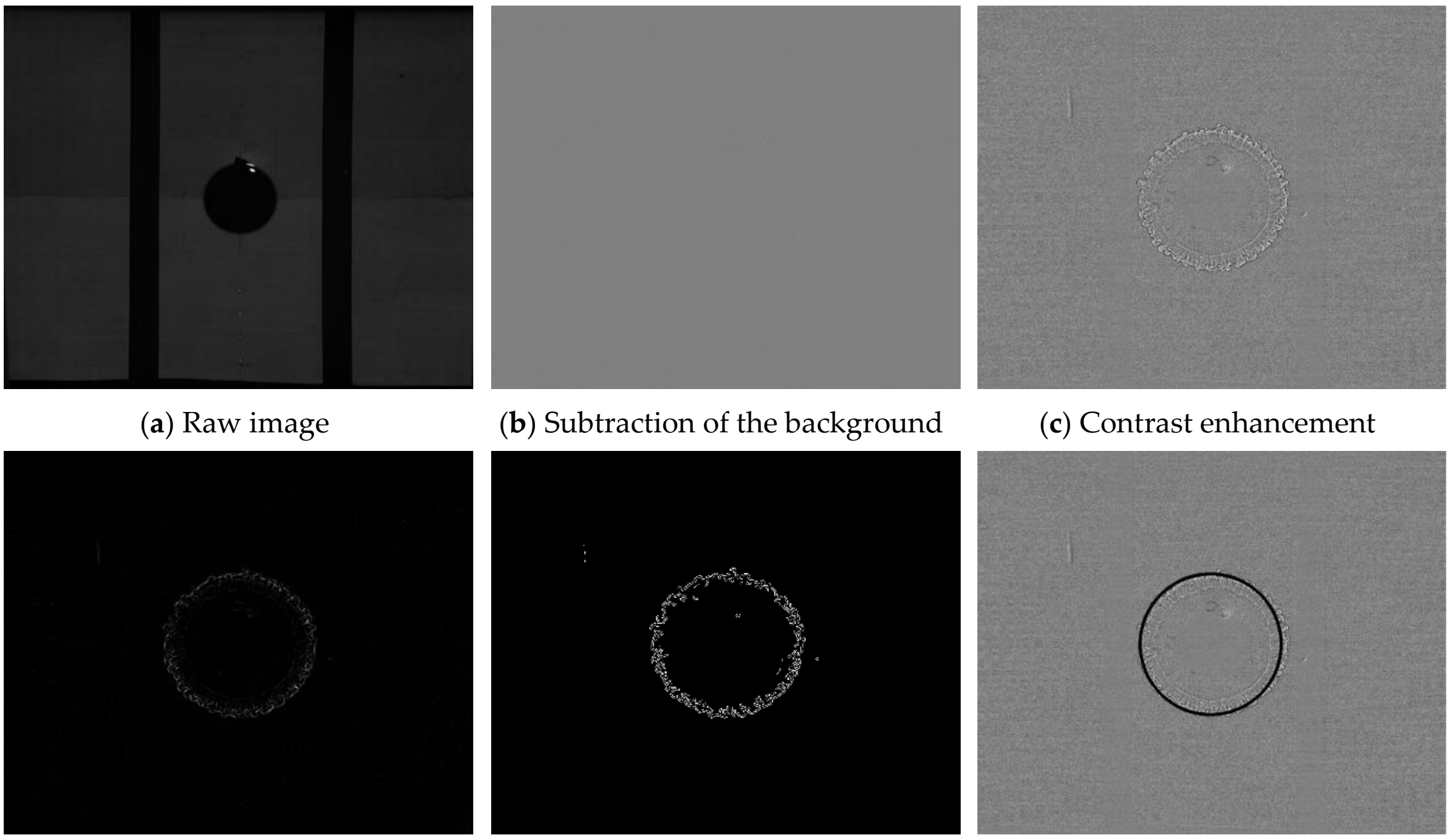

(**a**) Raw image (**b**) Subtraction of the background (**c**) Contrast enhancement

(**d**) Contour (first pass) (**e**) Contour (second pass) (**f**) Determination of the best-fit circle

**Figure 4.** Processing steps for the top-view shadowgraphs.

In the front view images, the cloud has a flat shape. Consequently, it is not possible to use strictly the same method to determine the radius of the cloud and the processing technique has been adapted. Instead of processing the whole image, only a horizontal line of pixels is processed. This line was taken five pixels above the floor which approximately corresponds to the height at which the head of the cloud reaches its maximum radius. The processing steps are illustrated in Figure 5. First, on each picture of the video, the line of pixels is extracted, and the background is subtracted. These lines are then assembled chronologically from top to bottom to form a kymograph (Figure 5a). The contrast of the kymograph is then enhanced using the CLAHE method (Figure 5b). The same two-edge filters [38,39] are then applied to get a binary image of the contours of the cloud. As artefacts of the background remained materialized as vertical streaks, only the vertical gradient of the image was kept (Figure 5c). In the last step, the left and right boundaries of the clouds are determined (Figure 5d). For that, the boundaries are first initialized as the boundaries of the container. Then the kymograph is processed line by line from top to bottom. The position of the left boundary on one line is chosen as the one of the leftmost contours detected on this line, in the limit of 50 pixels around the position calculated at the previous line. The symmetrical operation is performed for the right boundary. As for the processing of the images of the top view, the extraction of the contours is limited by the initial contrast of the images. Consequently, it is increasingly difficult to determine the boundaries of the cloud as it spreads. This is particularly visible in Figure 5c as fewer contours are detected as time evolves.

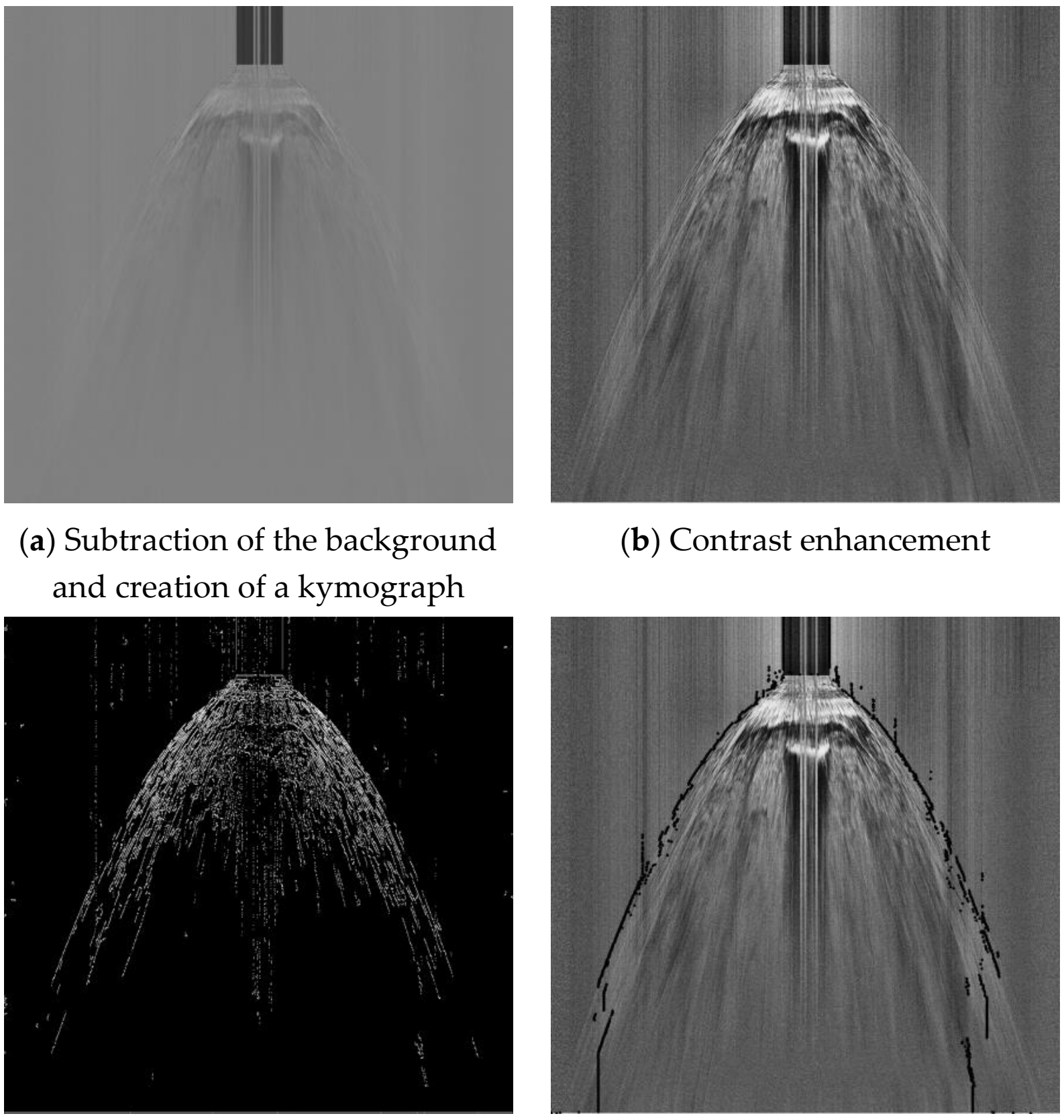

(**a**) Subtraction of the background and creation of a kymograph

(**b**) Contrast enhancement

(**c**) Contour

(**d**) Determination of the cloud boundaries

**Figure 5.** Processing steps for the front-view shadowgraphs.

### *4.4. Experimental Configurations and Data*

For all the tested configurations, the atmospheric pressure $p_a$ and temperature $T_a$ were measured near the vessel at the time of the trial. It allows calculating the actual density of krypton by using the ideal gas law:

$$\rho_{Kr} = \frac{p_a}{T_a}\frac{W_{Kr}}{\hat{R}}, \tag{16}$$

where $W_{Kr} = 83.8 \cdot 10^{-3}$ kg mol$^{-1}$ is the molar mass of krypton and $\hat{R} = 8.314$ J/mol/K is the ideal gas constant. Pressure and temperature equilibria between krypton and atmospheric air are assumed. It could be objected that, since krypton is expanded from a high-pressure bottle, it could be colder than the atmosphere. However, krypton flows slowly in a several-meter long pipe before arriving in the container and this pipe acts as a heat exchanger. Equation (16) is only used to provide the order of magnitude of the cloud Froude and Reynolds numbers for the various container heights as will be seen in Section 4.7. Consequently, assuming pressure and temperature equilibria between krypton and air seems reasonable.

Given the duration of the filling process, $t_{Kr}$, and the krypton volumetric flow rate, $Q_{Kr}$, it is possible to determine the actual volume of krypton inside the container and then the initial height of krypton:

$$V_{Kr} = Q_{Kr} t_{Kr}, \tag{17}$$

$$h_{Kr} = \frac{V_{Kr}}{\pi R_c^2}, \tag{18}$$

with $R_c$ the internal radius of the container. The last expression supposes that krypton does not mix substantially with air during the filling process. This point has been checked thanks to the shadowgraphs taken at the moments when the container starts to fall. Their analysis indeed shows that the height at which the separation between krypton and air begins to be visible corresponds closely to the calculated value of the height of krypton. Consequently, even if krypton happens to mix with air, the mixing is not significant.

All the experimental configurations are listed in Table 1. Three initial heights of krypton were tested and at least six trials were carried out for each height for repeatability. In the table and in the following, the experimental cases are referenced using a notation of the form "hX" where "h" is the target value for the height of krypton in centimeter and "X" is an alphabetical index. The column $\Delta h_{Kr}$ corresponds to the difference between the target value of the initial height of krypton and the actual value. The uncertainties indicated with the different parameters are calculated by assuming that the duration of the filling process and the krypton flow rate are known respectively with a precision of 1 s and 0.01 L min$^{-1}$, which corresponds to the precision of the flowmeter. For every case, the actual initial height of krypton is close to the target value and the difference between the two is always below the calculated uncertainty, which confirms repeatability of the initial conditions.

**Table 1.** Summary of the experimental configurations and associated parameters.

| **Case** | $t_{Kr}$ | $Q_{Kr}$ **[L·min$^{-1}$]** | $T_a$ **[°C]** | $p_a$ **[mbar]** | $\rho_{Kr}$ **[kg·m$^{-3}$]** | $V_{Kr}$ **[L]** | $h_{Kr}$ **[cm]** | $\Delta h_{Kr}$ **[cm]** |
|---|---|---|---|---|---|---|---|---|
| 5A | 1 min 34 s | 0.25 ± 0.01 | 16.8 | 976 | 3.39 | 0.39 ± 0.02 | 5.0 ± 0.3 | 0.0 ± 0.3 |
| 5B | 1 min 35 s | 0.24 ± 0.01 | 16.9 | 976 | 3.39 | 0.38 ± 0.02 | 4.9 ± 0.3 | −0.1 ± 0.3 |
| 5C | 1 min 34 s | 0.25 ± 0.01 | 17.8 | 981 | 3.40 | 0.39 ± 0.02 | 5.0 ± 0.3 | 0.0 ± 0.3 |
| 5D | 1 min 35 s | 0.24 ± 0.01 | 17.7 | 981 | 3.40 | 0.38 ± 0.02 | 4.9 ± 0.3 | −0.1 ± 0.3 |
| 5E | 1 min 35 s | 0.24 ± 0.01 | 17.6 | 981 | 3.40 | 0.38 ± 0.02 | 4.9 ± 0.3 | −0.1 ± 0.3 |
| 5F | 1 min 36 s | 0.25 ± 0.01 | 18.1 | 981 | 3.39 | 0.40 ± 0.02 | 5.1 ± 0.3 | 0.1 ± 0.3 |
| 5G | 1 min 34 s | 0.24 ± 0.01 | 18.0 | 981 | 3.40 | 0.38 ± 0.02 | 4.9 ± 0.3 | −0.1 ± 0.3 |
| 5H | 1 min 35 s | 0.24 ± 0.01 | 18.3 | 980 | 3.39 | 0.38 ± 0.02 | 4.9 ± 0.3 | −0.1 ± 0.3 |
| 10A | 3 min 9 s | 0.25 ± 0.01 | 18.0 | 980 | 3.39 | 0.79 ± 0.03 | 10.1 ± 0.4 | 0.1 ± 0.4 |
| 10B | 3 min 9 s | 0.24 ± 0.01 | 17.9 | 980 | 3.39 | 0.76 ± 0.03 | 9.8 ± 0.4 | −0.2 ± 0.4 |
| 10C | 3 min 8 s | 0.24 ± 0.01 | 15.9 | 978 | 3.41 | 0.75 ± 0.03 | 9.6 ± 0.4 | −0.4 ± 0.4 |
| 10D | 3 min 8 s | 0.25 ± 0.01 | 16.5 | 978 | 3.40 | 0.78 ± 0.03 | 10.0 ± 0.4 | 0.0 ± 0.4 |
| 10E | 3 min 9 s | 0.24 ± 0.01 | 16.3 | 978 | 3.41 | 0.76 ± 0.03 | 9.8 ± 0.4 | −0.2 ± 0.4 |
| 10F | 3 min 8 s | 0.25 ± 0.01 | 16.7 | 978 | 3.40 | 0.78 ± 0.03 | 10.0 ± 0.4 | 0.0 ± 0.4 |
| 15A | 4 min 46 s | 0.24 ± 0.01 | 15.0 | 975 | 3.41 | 1.14 ± 0.05 | 14.6 ± 0.6 | −0.4 ± 0.6 |
| 15B | 4 min 45 s | 0.24 ± 0.01 | 15.1 | 975 | 3.41 | 1.14 ± 0.05 | 14.6 ± 0.6 | −0.4 ± 0.6 |
| 15C | 4 min 44 s | 0.25 ± 0.01 | 15.4 | 975 | 3.41 | 1.18 ± 0.05 | 15.1 ± 0.6 | 0.1 ± 0.6 |
| 15D | 4 min 44 s | 0.25 ± 0.01 | 15.5 | 975 | 3.41 | 1.18 ± 0.05 | 15.1 ± 0.6 | 0.1 ± 0.6 |
| 15E | 4 min 44 s | 0.25 ± 0.01 | 16.3 | 975 | 3.40 | 1.18 ± 0.05 | 15.1 ± 0.6 | 0.1 ± 0.6 |
| 15F | 4 min 44 s | 0.25 ± 0.01 | 16.2 | 975 | 3.40 | 1.18 ± 0.05 | 15.1 ± 0.6 | 0.1 ± 0.6 |

### 4.5. Scaling Law

In order to compare the results obtained with the three container heights in the same frame, a scaling law will be used. The following parametrization [40] is used for the cloud radius $R$ and time $t$:

$$R \rightarrow \lambda = \frac{R}{\left(R_c^2 h_0\right)^{\frac{1}{3}}} = S_r R, \tag{19}$$

$$t \rightarrow \tau = \left(S_r g\left(\frac{\rho_0}{\rho_a} - 1\right)\right)^{\frac{1}{2}} t = S_t t, \tag{20}$$

where $\lambda$ and $\tau$ are the dimensionless scaled radius and time, $R_c$ is the internal radius of the container, $h_0$ is the initial height of the heavy gas cloud, $\rho_0$ is its initial density, $\rho_a$ is the atmospheric density and $g$ is the acceleration due to gravity.

Therefore, the cloud area $A$ can also be scaled according to the following relation:

$$A \rightarrow \Sigma = S_r^2 A. \tag{21}$$

### 4.6. Expected Evolution of the Cloud

In the present cylindrical experimental setup, the heavy gas is released almost instantaneously, and the subsequent flow is axisymmetric. Box models [1,2] have been developed in this specific context. Those are however irrelevant for more complex flow configurations. They are valid only in cylindrical and axisymmetric situation such as the present one. It is consequently worth introducing such models. The experimental time evolution of the cloud, in terms of disk-shaped area, is then to be compared with these specific models. In the simplest version, the heavy gas cloud is represented as a cylindrical volume of homogeneous and constant density gas spreading on the ground. The cloud volume is thus assumed to be constant. The expansion velocity is then given by:

$$\frac{dR}{dt} = Fr\sqrt{g\frac{\rho_{Kr} - \rho_a}{\rho_a} h_{Kr}}, \tag{22}$$

where $R$ and $h_{Kr}$ are the instantaneous radius and height of the cloud, $g$ is the acceleration due to gravity, $\rho_{Kr}$ and $\rho_a$ are the respective density of the heavy gas and air and $Fr$ is the cloud Froude number which represents the ratio between its kinetic energy and its gravitational potential energy. From this expression, it is possible to derive an evolution equation for the disk-shaped area of the cloud:

$$\frac{dA}{dt} = 2\pi Fr\sqrt{\frac{gV_{Kr}}{\pi}\frac{\rho_{Kr} - \rho_a}{\rho_a}}, \tag{23}$$

$$A(t) = 2\pi Fr\sqrt{\frac{gV_{Kr}}{\pi}\frac{\rho_{Kr} - \rho_a}{\rho_a}}t + A_0, \tag{24}$$

with $V_{Kr}$ being the (constant) cloud volume and $A_0$ the initial cloud area.

Seen from above, the present experimental cloud is expected a have a circular shape. According to Equation (24), its disk-shaped area $A(t)$ is to evolve linearly with time. Its radius $R(t)$ is consequently expected to present a parabolic time evolution. The aera and radius are indeed related through $A(t) = \pi R^2(t)$. As will be seen further in Sections 4.7 and 5.3, present experimental data verify indeed these behaviors.

In the right side of Equation (24), the Froude number is the only unknown. In the left side, the area of the cloud is known with the help of the cloud radius, experimentally measured by the methods (top and front views) presented previously and evolves linearly with time. The Froude number can consequently be extracted from Equation (24) by performing a linear regression.

Experimental results are presented in the following section. The corresponding Froude numbers are also provided. Let us recall that such simple model is only valid for cylindrical and symmetrical flows. It is here used for the sole purpose of checking the consistency of the experimental data. Indeed, such model cannot be used to address non-cylindrical flows or flows in complex geometries. The present two-layer shallow water model (1) is obviously more complex than Equation (24) but can deal with arbitrary geometrical domains. Its accuracy will be examined in Section 5.3.

### *4.7. Results and Discussion*

Shadowgraphs taken during trial 15E and showing the expansion of the krypton cloud are shown in Figure 6. These shadowgraphs have been processed using the method described in Section 4.3 in order to remove the background and increase the contrasts. When the container reaches its lowest position (first row of Figure 6), the cloud begins to fall. The interface between air and krypton wrinkles due to the apparition of Kelvin-Helmholtz [41,42] instabilities. Then, krypton spreads on the floor. Seen from above, the cloud keeps an overall circular shape which confirms the good symmetry of the gas release. Unlike front view images may suggest, the thickness of the cloud is not uniform. Instead, it has a thick vertical "head" and is thin and flat elsewhere. Rayleigh-Taylor instabilities [43,44] can be seen on the "head" of the cloud and are caused by both the difference of density between krypton and air and deceleration of the cloud.

| t (ms) | Front View | Top View |
|---|---|---|
| 224 | 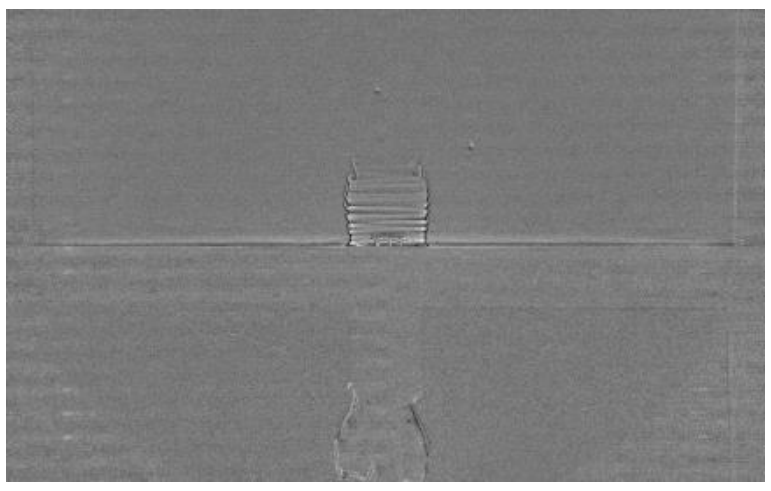 | 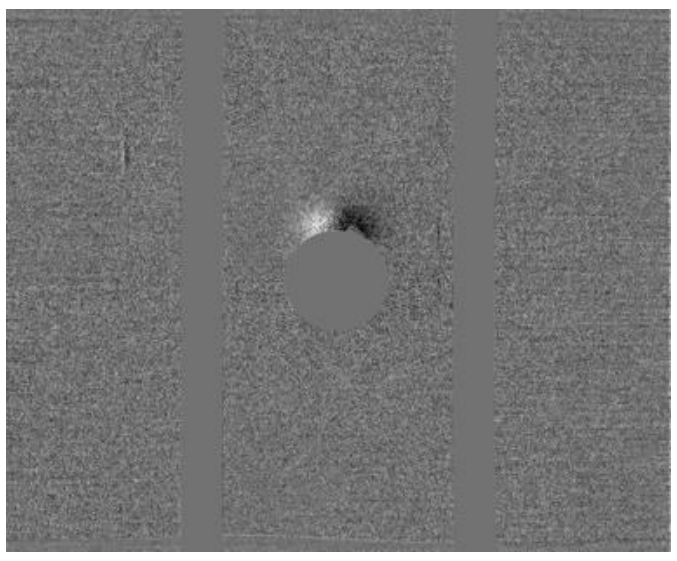 |
| 364 | 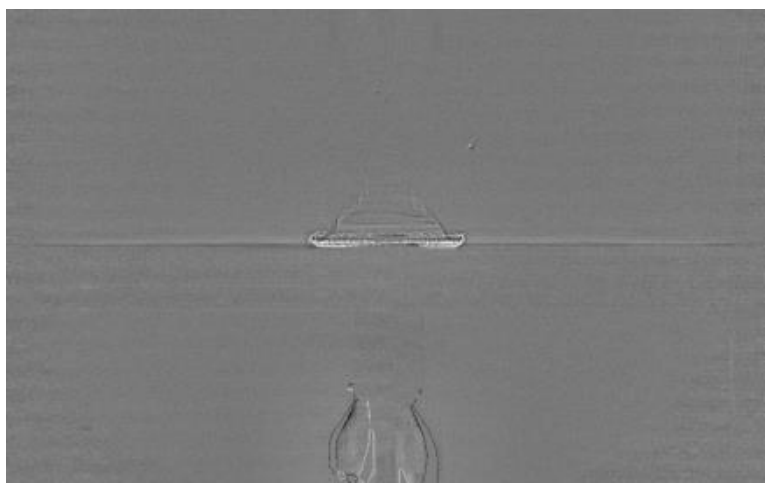 | 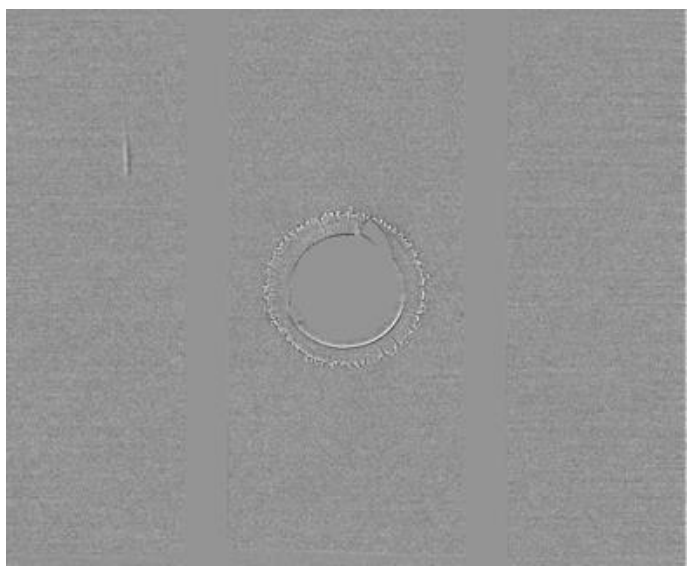 |
| 464 | 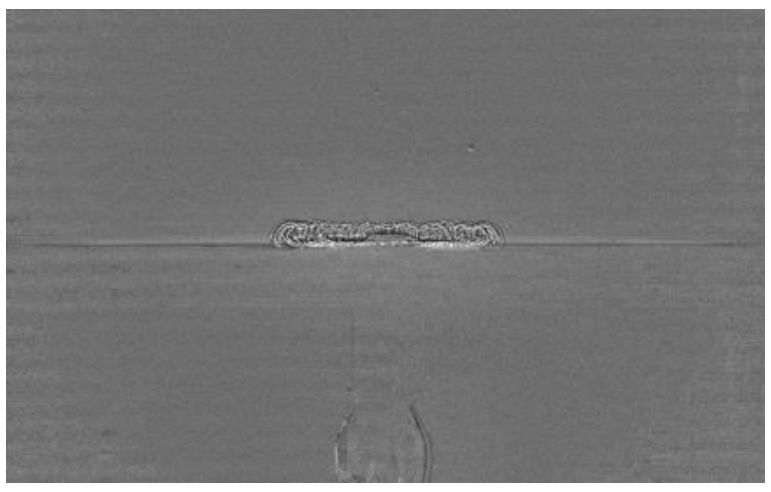 | 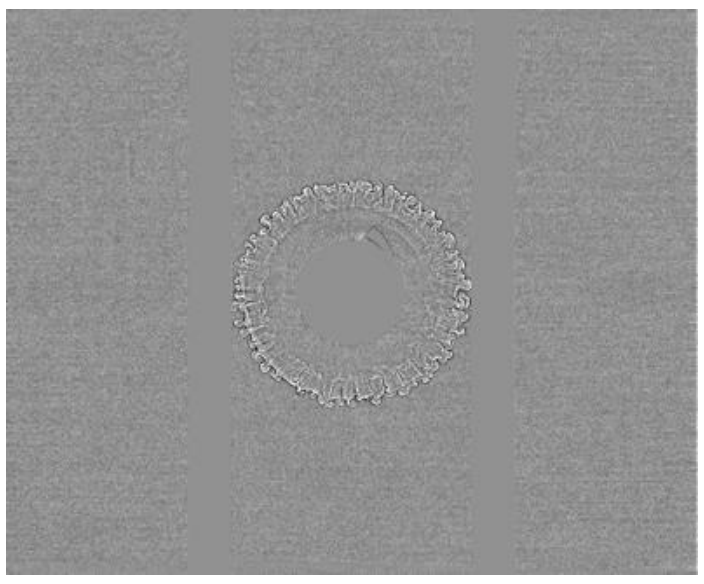 |
| 564 | 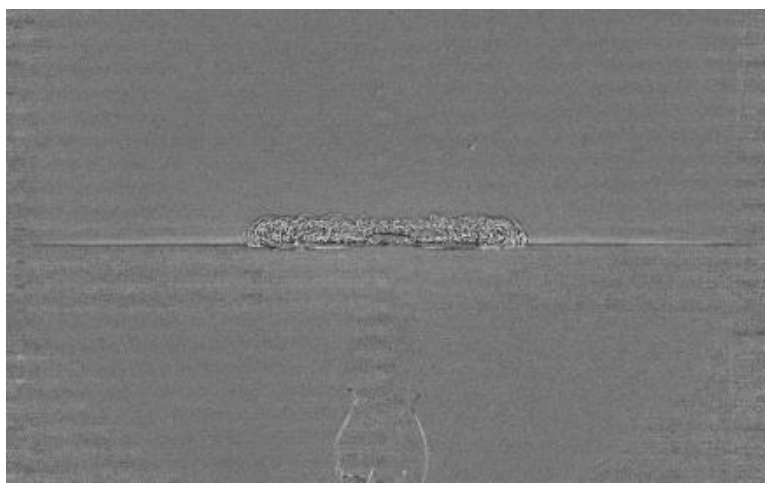 | 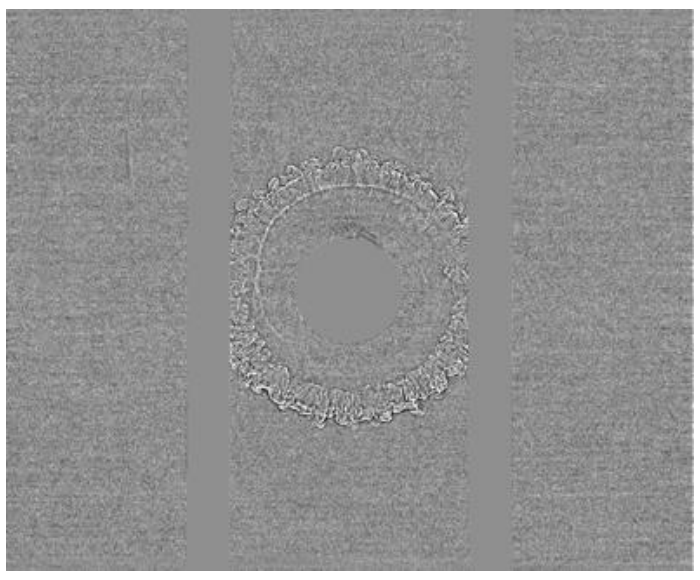 |

**Figure 6.** Shadowgraphs (front view and top view) of trial 15E showing the evolution of the krypton cloud. For each row, the front and top view pictures are taken at the same time. The origin of time corresponds to the beginning of the descent of the container.

Using the processing technique described in Section 4.3, the temporal evolution of the cloud radius has been extracted from the shadowgraphs and the evolution of its area has been deduced. The resulting data are plotted in Figure 7a for the top view and Figure 7b for the front view. For all sets of data, the origin of time has been manually corrected in order to correspond to the end of the descent of the container. As it can be seen, the experimental setup proved to be repeatable for each initial height of gas despite its apparent simplicity. Processing of the front view images leads to more scattered data but allows following the cloud for a longer time. It is worth noting that the top view image processing tends to underestimate the cloud area above 0.08 m$^2$. A slight inflection can for example be seen at 0.23 s on the graphs in Figure 7 corresponding to the initial height of 15 cm. Such inflection appears because part of the cloud is then located over the two beams that support the floor of the vessel, which hinders the determination of the circle.

The growth of the cloud's area is linear as expected (see Section 4.6). The slope of the curves however still depends on the Froude number. When the scaling law from Section 4.5 is applied, all the plots almost collapse on a single line as shown in Figure 8. These results confirm the good behavior of experimental data and their validity for the validation of the two-layer model (1).

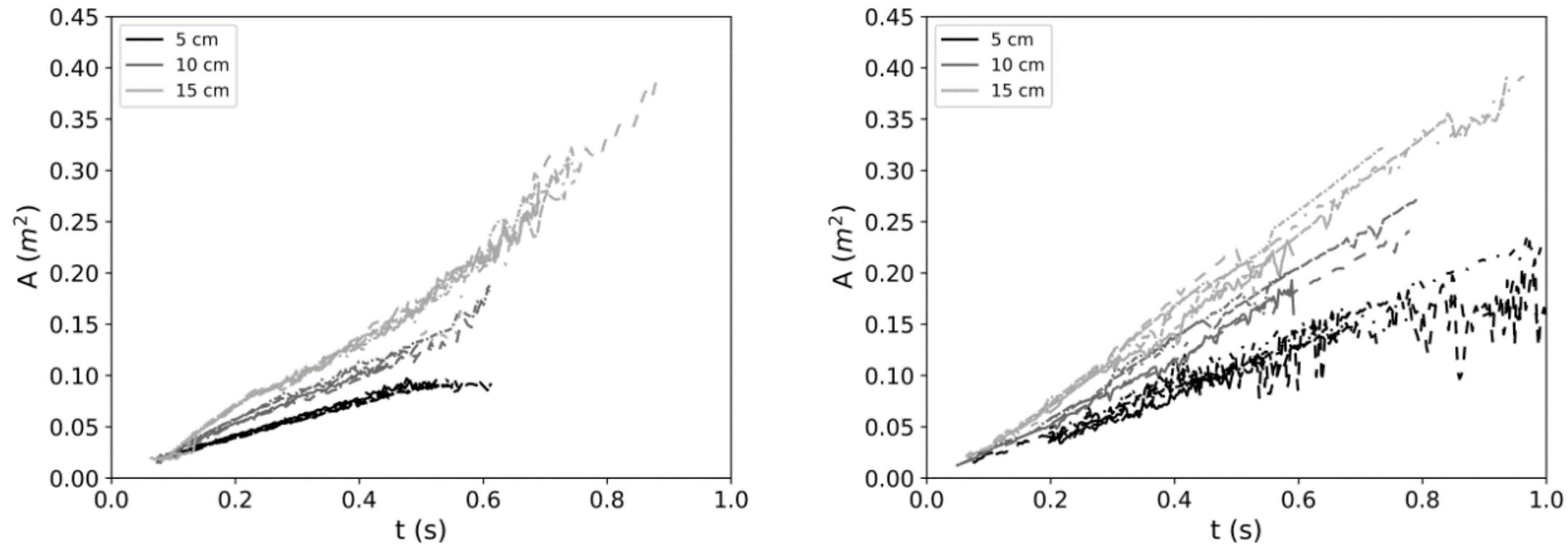

(**a**) Top view (**b**) Front view

**Figure 7.** Temporal evolution of the cloud area measured on the shadowgraphs.

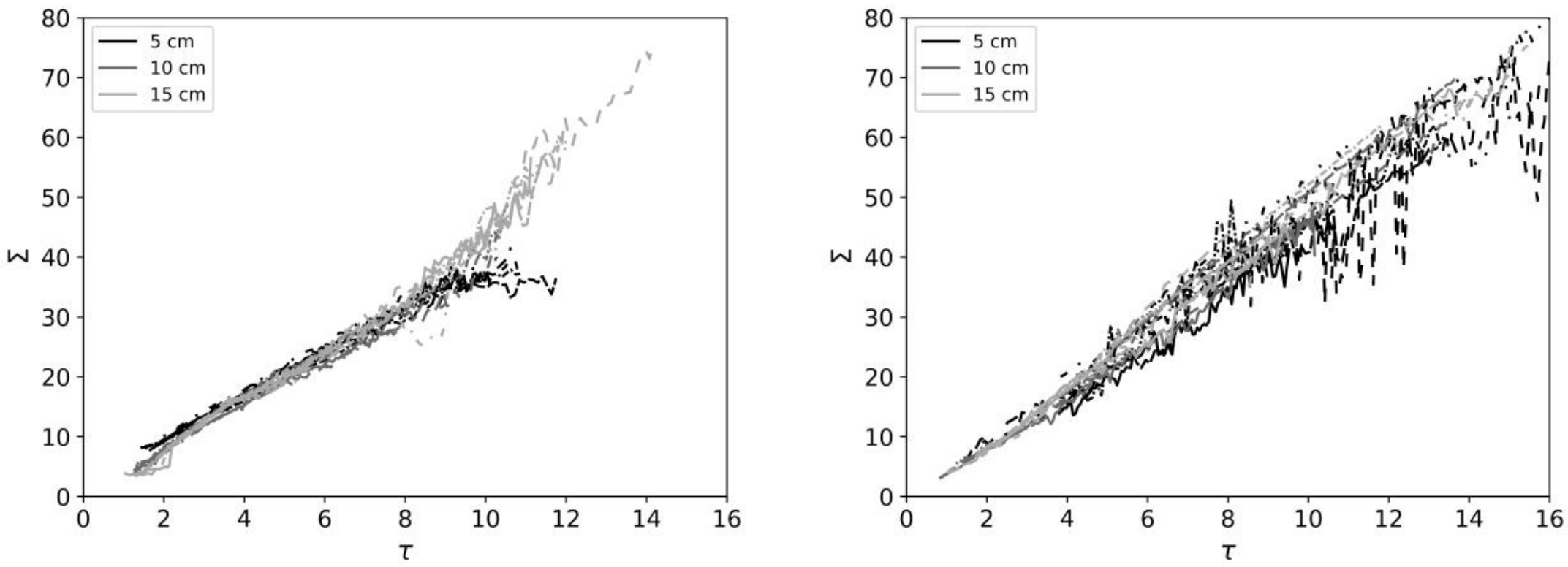

(**a**) Top view (**b**) Front view

**Figure 8.** Temporal evolution of the cloud area (defined by Equation (21)) measured on the shadowgraphs in scaled coordinates.

From the evolution of the cloud area, the Froude number has been determined for every configuration by performing a linear regression using Equation (24). Regarding front view images, as they could contain some noise, a filtering step was performed before applying linear regression. By comparing the determined cloud radius with the front view pictures, it was noticed that the noise corresponded mostly to an underestimation of the radius. This is due to turbulent structures located inside the cloud instead of its boundaries. Consequently, instead of performing directly linear regression on the data output by the processing script, it was decided to perform it on the upper envelope of the curves. This method appeared to be effective even on very noisy signals. The determined value of the Froude number is given for each configuration in Table 2. Top view data for trial 5B and trial 10B do not appear in the table because of data transfer issues.

**Table 2.** Froude number calculated from the processing of the shadowgraphs.

| Case | Froude Number [–] | | | |
|---|---|---|---|---|
| | Top View | Front View–Left | Front View–Right | Average |
| 5A | 0.61 | 0.80 | 0.70 | 0.70 |
| 5B | - | 0.79 | 0.65 | 0.72 |
| 5C | 0.63 | 0.72 | 0.78 | 0.71 |
| 5D | 0.58 | 0.80 | 0.63 | 0.67 |
| 5E | 0.61 | 0.82 | 0.75 | 0.73 |
| 5F | 0.61 | 0.71 | 0.58 | 0.63 |
| 5G | 0.58 | 0.71 | 0.69 | 0.66 |
| 5H | 0.61 | 0.76 | 0.73 | 0.70 |
| 10A | 0.63 | 0.78 | 0.81 | 0.74 |
| 10C | 0.57 | 0.79 | 0.73 | 0.70 |
| 10D | 0.62 | 0.77 | 0.81 | 0.73 |
| 10E | 0.61 | 0.86 | 0.93 | 0.80 |
| 10F | 0.63 | 0.86 | 0.71 | 0.73 |
| 15A | 0.80 | 0.80 | 1.01 | 0.87 |
| 15B | 0.81 | 1.00 | 0.84 | 0.88 |
| 15C | 0.77 | 0.81 | 0.93 | 0.84 |
| 15D | 0.81 | 0.80 | 0.85 | 0.82 |
| 15E | 0.72 | 0.90 | 0.82 | 0.81 |
| 15F | 0.75 | 0.87 | 0.96 | 0.86 |

From the data given in Table 2, the average Froude number was calculated for each initial height of krypton. The results are given in Table 3. The associated standard deviation was also calculated to give an estimation of the uncertainty. As it can be noticed, the average Froude number tends to increase with the initial cloud height, which shows that the conversion of potential energy into kinetic energy tends to be more efficient when the cloud is initially higher.

**Table 3.** Average Froude number for each initial height and associated standard deviation.

| $h_{Kr}$ [cm] | Average $Fr$ [–] | Standard Deviation [–] |
|---|---|---|
| 5 | 0.69 | 0.03 |
| 10 | 0.74 | 0.03 |
| 15 | 0.85 | 0.03 |

As the configuration of the current trials is similar to the one studied at Thorney Island [35] and only differs in scale, it is of interest to compare the results of both experiments. The trials of Thorney Island (Phase I) consisted of the instantaneous release of a dense gaseous mixture of dichlorofluoromethane and nitrogen initially contained in a twelve-sided, 14-m wide and 13-m high container, which corresponds to an initial volume of 2000 m$^3$. The evolution of the cloud radius and area was estimated thanks to aerial photographs taken from a helicopter located 300 m above the test site [35]. The analysis of the collected data allowed to determine that the cloud Froude number was 1.05 ± 0.12 [1]. This value is higher than the ones found during the current trials, which tends to corroborate the previous observation made about the Froude number.

It is highly likely that viscous effects play a more prominent role at small scale than at full scale. The cloud Reynolds number, representing the ratio between inertial and viscous effects, reads,

$$\mathrm{Re} = \frac{\rho_a u R}{\mu_a}, \tag{25}$$

where $u$ is the cloud front velocity and $\mu_a$ is the viscosity of air. According to Equation (22), the cloud's expansion velocity obeys,

$$R\frac{dR}{dt} = Ru = Fr\sqrt{\frac{gV_{Kr}}{\pi}\frac{\rho_{Kr}-\rho_a}{\rho_a}} = K. \quad (26)$$

Hence, the expression of the cloud Reynolds number becomes,

$$Re = \frac{\rho_a}{\mu_a}K, \quad (27)$$

with $\rho_a = \frac{p_a}{T_a}\frac{W_a}{\hat{R}}$ and $W_a = 29\cdot 10^{-3}$ kg.mol$^{-1}$. For the current trials, the value of the constant is $K \simeq 0.06$ m$^2$/s for $h_{Kr} = 15$ cm, $K \simeq 0.05$ m$^2$/s for $h_{Kr} = 10$ cm and $K \simeq 0.03$ m$^2$/s for $h_{Kr} = 5$ cm. The cloud Reynolds number is thus smaller for small heights of container and viscous effects are indeed more prominent in that case. In the case of the Thorney Island trials [35], $K \simeq 68$ m$^2$/s, which tends to confirm the previous conclusion.

In the frame of the validation of the two-layer model, it means that, although the current experimental results can indeed be used as a reference for comparison, attention shall be paid in future work to scale effects.

Numerical results of the present two-layer shallow water model (1) are compared to experimental data in the next section. In this specific cylindrical symmetry context, simpler models such as (23) and (24), provide good results, assuming known Froude number. However, such simplified models become irrelevant for non-cylindrical situations and are unable to deal with more complex geometries. The two-layer shallow water approach becomes more appropriate.

## 5. Numerical Results with Two-Layer Model and Laboratory Experiments Comparison

### *5.1. Two-Layer Flow Model with Cylindrical Symmetry*

The present laboratory experiments involve a cylindrical gas column. Consequently, cylindrical symmetry is introduced in the 1D equations of System (1). The corresponding system is presented in Appendix A.

### *5.2. Velocity Relaxation*

Numerical experiments have been carried out in order to find the best values of the interfacial area $A_I$ (see Section 2, Equation (5)) for each dam-break configuration. For each test, multiple values of $A_I$ (constant during time step $\Delta t$ ) have been necessary to match computational and experimental results. Based on these numerical experiments, the following function has been determined:

$$A_I(t) = A_I(h_0)e^{-\frac{t}{\tau}} + A_{I0}. \quad (28)$$

$A_I(h_0)$ represents the initial interfacial area. It consequently depends on the initial height $h_0$ of the krypton column. $\tau$ represents the "system's time constant" and $A_{I0}$ is a positive constant ensuring that the interfacial area $A_I(t)$ remains non-zero during time evolution. The present drag function may be conveniently written with the following parameters,

$$A_I\left(t\right) = a.h_0 e^{-b.t} + c, \quad (29)$$

where $h_0$ is the initial height of the krypton column and $a, b, c$ are three positive constants. The interfacial area is then expressed under the form of a correlation depending on time and requiring three parameters.

According to numerical experiments, parameters $a = 0.2\ m^{-1}$, $b = 12.5\ s^{-1}$, $c = 5.10^{-5}$ are a fair choice. The present interfacial law consequently reads,

$$A_I(t) = 0.2\ h_0 e^{-12.5\ t} + 5.10^{-5}, \tag{30}$$

for all initial conditions of two-layer dam-break test problems.

*5.3. Results and Comparisons*

The experiments consist of a cylindrical krypton column of diameter 10 cm. Under gravity effects, gas dispersal occurs into the surrounding atmospheric air. Experimental results are provided in terms of position (radius $R$ ) of the krypton front. Three initial column heights are used, namely 5, 10 and 15 cm. The experimental apparatus consists of a glass-made cubic enclosure of dimension $1.2 \times 1.2 \times 1.2$ m³. However, for symmetry reasons, the following 1D computations only consider half of the domain. These initial data correspond to the ones of the experiments of Section 4.

Initial conditions consequently consist of a 0.6-m wide computational domain with a discontinuity located at $x_0 = 0.05$ m. On the left of the discontinuity, initial heights are $h_1 = 0.05,\ 0.10,\ 0.15$ m and $h_2 = 1.2 - h_1$. On the right, the initial heights are $h_1 = 10^{-6}$ m and $h_2 = 1.2$ m. Initial densities are $\rho_1 = 3.506$ kg/m³ for the krypton and $\rho_2 = 1.29$ kg/m³ for the air. The gravity constant is $g = 9.81$ m/s². Fluids are initially at rest $u_1 = u_2 = 0$ m/s. Reflective boundary conditions are used (symmetry condition on the left and wall on the right). Acoustic impedances $Z_k = \rho_k c_k = (\rho_k c_k)_{atm}$ are considered constant in the drag force modeling and are set to their atmospheric values: $\rho_1 = 3.506$ kg/m³ and $c_1 = 218$ m/s for the krypton and $\rho_2 = 1.29$ kg/m³ and $c_2 = 340$ m/s for the air. The atmospheric pressure is $p_0 = 101325$ Pa. The "sound speed" numerical parameters are $\theta_1 = \theta_2 = 2$ for both fluids during the hyperbolic step.

Computed results' accuracy is improved with the Monotonic-Upwind-Scheme-for-Conservation-Laws (MUSCL) type second order scheme (see for example Chiapolino et al., 2017 [27]) using van Leer's limiter [45]. The computational domain is made of 1000 cells $(M = 1000)$. The CFL number is $CFL = 0.8$. The computed position of the krypton front is specifically determined according to,

$$h_{1,i} h_{2,i} > \varepsilon. \tag{31}$$

This indicator notifies height variations at cell center $i$. Relation (31) is used from right $(i = M)$ to left $(i = 1)$. When a height variation of $\varepsilon = 10^{-10}$ m² is recorded, the krypton front position is considered at $\frac{x_i + x_{i+1}}{2}$ where $x_i$ and $x_{i+1}$ represent locations of the cell centers $i$ and $i+1$.

Numerical results are compared to experimental data in Figures 9–11. Experimental results are extracted from multiple series of tests resulting in various (grey lines) curves in the following graphs. The dam-break test problem is first considered with a 5-cm high initial krypton column. Corresponding results are shown in Figure 9.

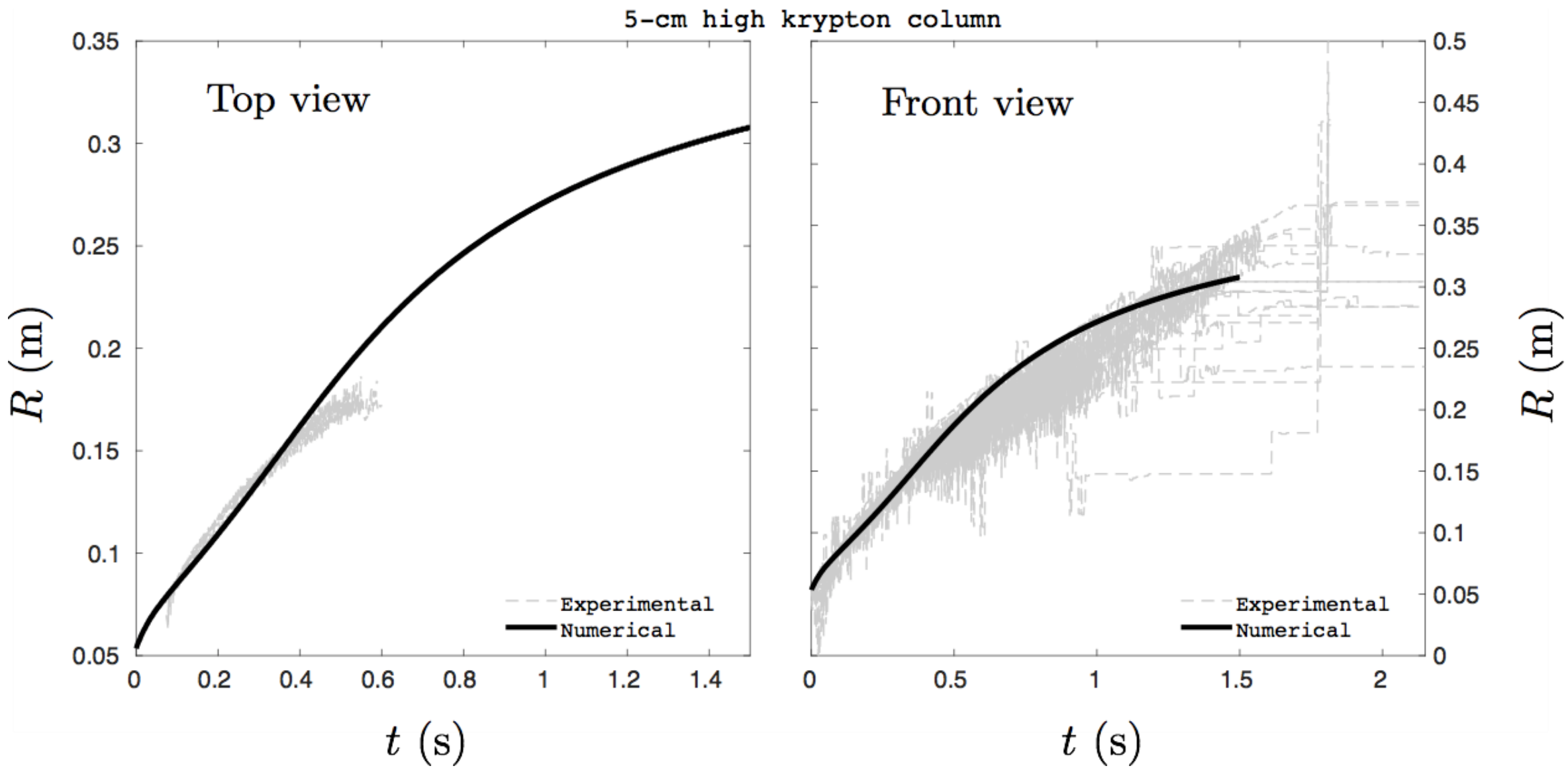

**Figure 9.** Temporal evolution of the krypton front position. Comparison of the one-dimensional solution of the present two-layer shallow water model versus experimental data. The initial krypton column is 5-cm high. Computed results are in good agreement with experimental measurements.

Good agreement between numerical and experimental results is observed. The test is now repeated with a 10-cm high initial krypton column. Corresponding results are shown in Figure 10.

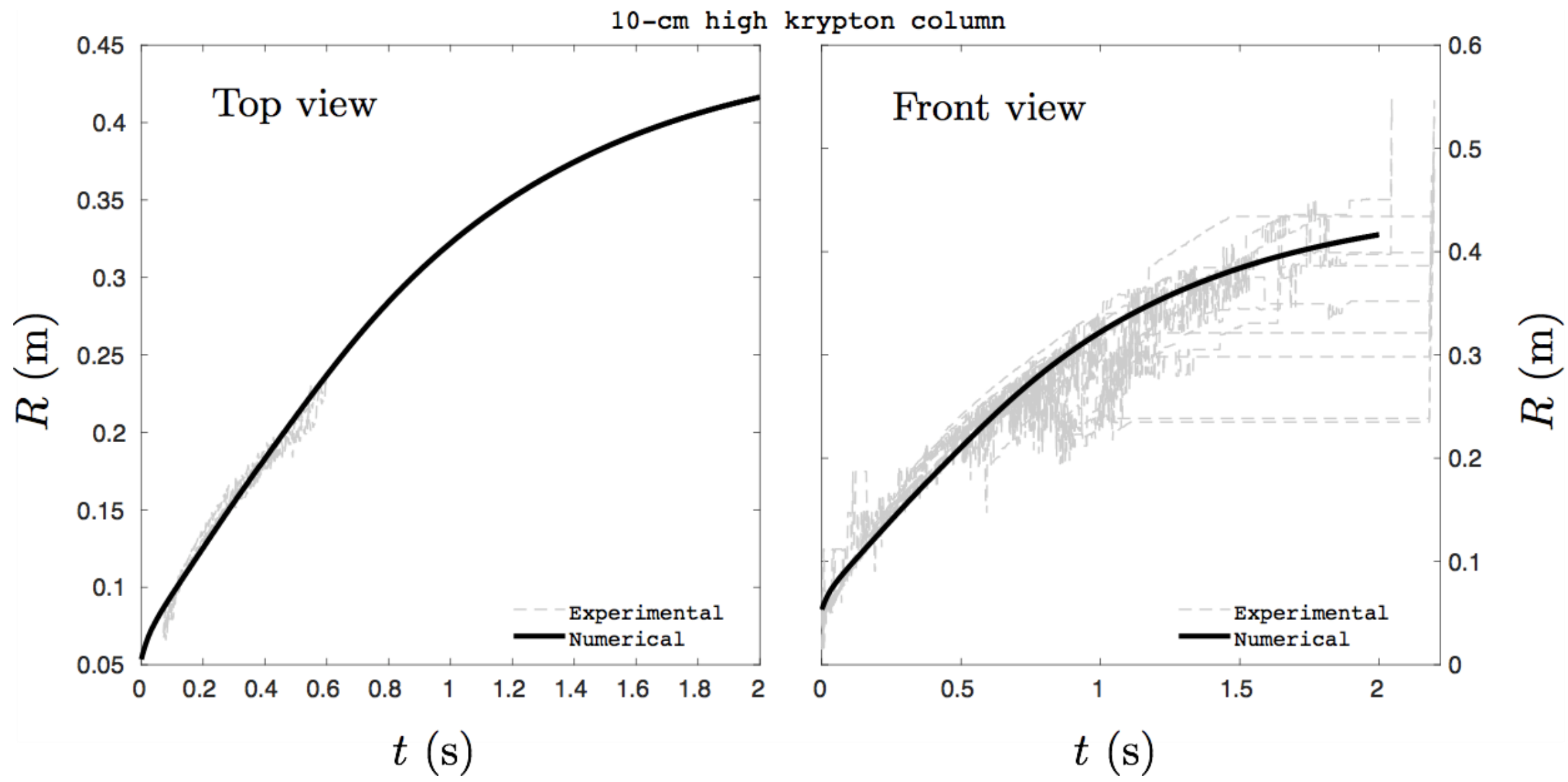

**Figure 10.** Temporal evolution of the krypton front position. Comparison of the one-dimensional solution of the present two-layer shallow water model versus experimental data. The initial krypton column is 10-cm high. Computed results are in good agreement with experimental measurements.

Again, good agreement is observed. The test is finally repeated with a 15-cm high initial krypton column. Corresponding results are shown in Figure 11.

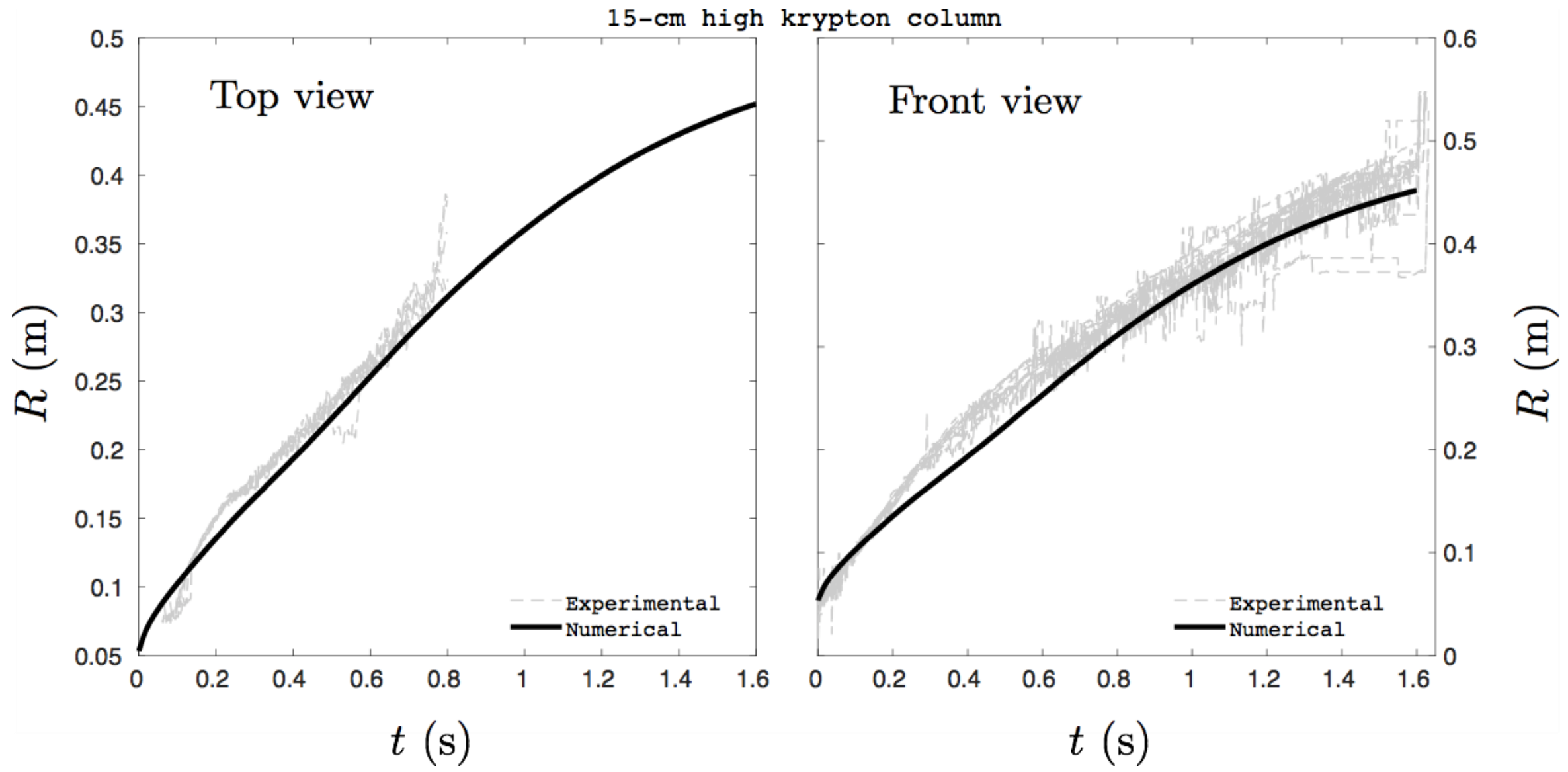

**Figure 11.** Temporal evolution of the krypton front position. Comparison of the one-dimensional solution of the present two-layer shallow water model versus experimental data. The initial krypton column is 15-cm high. Computed results are in good agreement with experimental measurements.

Once more, good agreement between numerical and experimental results is obtained. Those results indicate that the present hyperbolic two-layer shallow water model accurately reproduces experimental results at the price of few parameters related to the drag function.

In the three studied configurations, Relation (30) inserted in the acoustic drag force model (5) enables fairly accurate comparison of computed and experimental results. These small-scale comparisons show model's ability to predict large-scale dispersal of dense gases when drag effects are properly adjusted.

### 5.4. Concluding Remarks

The present approach based on two-layer shallow water equations appears promising. The agreement between experimental and computational results is good. It is obtained nowadays at the price of a few parameters related to the interfacial aera function, summarizing drag and turbulence effects at the interface separating the dense gas and the ambient one. A more sophisticated interfacial area will be needed in the future to address more complex situations and large-scale configurations. However, the present work shows that the two-layer approach has nice predicting abilities.

Topography effects, encountered in urban places or hilly grounds, have not been considered for the sake of simplicity. The model with its drag modeling has then been successfully compared to experimental data in its simplest form. This step was necessary before addressing more complex configurations.

Nevertheless, introduction of topography effects into two-layer shallow water equations is possible and will be part of future investigations. This flexibility is real asset of the two-layer shallow water approach, compared to other simplified box or SLAB-type models that are inappropriate for describing a gas cloud evolving in complex topographies [4].

The multi-D CFD approach is able to account for complex terrains but is much more complex and much more time-consuming than the present two-layer shallow water strategy. The corresponding equations have been recently revisited in Chiapolino and Saurel, 2018 [8] with an unconditionally hyperbolic formulation. In the same contribution, the importance of the two-layer property of the model, resulting in interaction terms, is illustrated. Indeed, in the present heavy-gas dispersal context, a single-layer model is too inaccurate, at least without extra terms [4]. Numerical resolution of the present two-layer model is about one million times faster than common multi-D CFD computations [8] and is robust and accurate if properly adjusted.

## 6. Conclusions

Heavy-gas dispersion modeling has been addressed in this paper. From a numerical point of view, one of the difficulties is related to long-time and large-scale computations, for instance a city scale, while providing accurate results with reasonable CPU time.

The present method differs from models commonly used in heavy-gas dispersion [7] as it relies on a two-layer shallow water approach. This approach seems more appropriate than box or SLAB-type models to address heavy-gas dispersal in urban places or hilly grounds. Indeed, among the aforementioned methods, introduction of topography effects is possible only with the two-layer approach. Such task has not been addressed in the present paper for the sake of simplicity but shall be part of future investigations. Computed results are indeed promising and highlight the predicting abilities of the two-layer shallow water approach.

Interactions between the ambient fluid and the dense layer is significant in this heavy-gas dispersal context and a single-layer model is inaccurate [8], at least without extra terms [4]. In the present contribution, these fluid interactions are accounted for through the non-conservative terms of the recent hyperbolic two-layer shallow water model of Chiapolino and Saurel, 2018 [8] and through drag effects that play a significant role. Particularly, creation of an interfacial aera is significant as it summarizes turbulent structures appearing at the interface and controls the dispersion speed of the dense gas cloud.

A set of laboratory experiments involving a dense gas has been specifically designed to address determination of the interfacial aera. Thanks to the experiments, an interfacial area function with few parameters has been determined and yields a good agreement between numerical and experimental results.

Computed results highlight the ability of the two-layer shallow water approach to address heavy-gas dispersal. However, a more sophisticated determination of the interfacial area will be needed in the future to address real configurations such as gas dispersal in a city, for instance. This interfacial area modeling is an important research work mainly related to friction between fluid layers that depends on interfacial instabilities. Such topics shall also be part of future investigations.

Most two-layer shallow water models are only conditionally hyperbolic, resulting in serious consequences both for wave propagation, which becomes ill-posed, and for the design of numerical methods. The present two-layer shallow water is unconditionally hyperbolic and consequently presents a solid base for many future developments. Indeed, in addition to the above-mentioned prospects, this research work can be continued in many other directions. Among them are the extension of the hyperbolic model to more than two fluids, the introduction of mass transfer between layers describing ascending effects or the modeling of thin and heterogeneous layers such as the free surface of the ocean (water-air mixing zone).

**Acknowledgments:** The authors wish to acknowledge the technical teams from CEA Gramat for their contribution to the conception of the experimental setup and the success of the experiments.

## Appendix A

*Appendix A.1. Two-Layer Flow Model with Cylindrical Symmetry*

In the present work, cylindrical symmetry is considered through geometric source terms resulting from the integration of the two-layer shallow water equations on a cylindrical control volume. For the sake of space restriction, calculation details are omitted and only the final system is given as

$$\begin{cases} \dfrac{\partial h_1}{\partial t} + u_1 \dfrac{\partial h_1}{\partial r} = 0, \\ \dfrac{\partial (h_1\rho_1)}{\partial t} + \dfrac{\partial (h_1\rho_1 u_1)}{\partial r} = -\dfrac{h_1\rho_1 u_1}{r}, \\ \dfrac{\partial (h_1\rho_1 u_1)}{\partial t} + \dfrac{\partial \left( h_1\rho_1 u_1^2 + h_1 (p_1 - p_0) + \frac{1}{2}\rho_1 g h_1^2 \right)}{\partial r} = \rho_2 g h_2 \dfrac{\partial h_1}{\partial r} - \dfrac{h_1\rho_1 u_1^2}{r}, \\ \dfrac{\partial h_2}{\partial t} + u_2 \dfrac{\partial h_2}{\partial r} = 0, \\ \dfrac{\partial (h_2\rho_2)}{\partial t} + \dfrac{\partial (h_2\rho_2 u_2)}{\partial r} = -\dfrac{h_2\rho_2 u_2}{r}, \\ \dfrac{\partial (h_2\rho_2 u_2)}{\partial t} + \dfrac{\partial \left( h_2\rho_2 u_2^2 + h_2 (p_2 - p_0) + \frac{1}{2}\rho_1 g h_1^2 \right)}{\partial r} = -\rho_2 g h_2 \dfrac{\partial h_1}{\partial r} - \dfrac{h_2\rho_2 u_2^2}{r}. \end{cases} \tag{A1}$$

Let us recall that index 1 denotes the heavy fluid (krypton) and index 2 the light one (air). Here, $r$ is the radial direction and the terms $-\dfrac{h_k\rho_k u_k}{r}$ and $-\dfrac{h_k\rho_k u_k^2}{r}$ are discretized, after the hyperbolic step, as:

$$\begin{aligned} -\frac{h_k\rho_k u_k}{r} &\to -\frac{1}{R_i}\frac{(h_k\rho_k u_k)^*_{i+1/2} + (h_k\rho_k u_k)^*_{i-1/2}}{2}, \\ -\frac{h_k\rho_k u_k^2}{r} &\to -\frac{1}{R_i}\frac{(h_k\rho_k u_k^2)^*_{i+1/2} + (h_k\rho_k u_k^2)^*_{i-1/2}}{2}, \end{aligned} \tag{A2}$$

where the radius $R_i$ is the location of the center of cell $i$ and the superscript $*$ denotes the Riemann problem solution (Section 3). Pressure relaxation is then achieved as in Section 2 followed by velocity relaxation. This step is addressed hereafter.

*Appendix A.2. Velocity Relaxation*

Drag effects are now included in the momentum equations (as seen in Section 2) as source terms:

$$\begin{cases} \dfrac{\partial(h_1\rho_1)}{\partial t} = 0, \\ \dfrac{\partial(h_1\rho_1 u_1)}{\partial t} = \dfrac{Z_1 Z_2}{Z_1 + Z_2} A_I (u_2 - u_1), \\ \dfrac{\partial(h_2\rho_2)}{\partial t} = 0 \\ \dfrac{\partial(h_2\rho_2 u_2)}{\partial t} = -\dfrac{Z_1 Z_2}{Z_1 + Z_2} A_I (u_2 - u_1). \end{cases} \tag{A3}$$

When the interfacial area $A_I$ is considered as constant, as well as acoustic impedances $Z_k = \rho_k c_k = (\rho_k c_k)_{atm}$ solution of System (A3) is determined explicitly as:

$$\begin{cases} (u_2)_t = \dfrac{(h_1\rho_1 u_1 + h_2\rho_2 u_2)^n + (h_1\rho_1)^n (u_2 - u_1)^n e^{(cst\ A_I\ \Delta t)}}{(h_1\rho_1)^n + (h_2\rho_2)^n}, \\ (u_1)_t = (u_2)_t - (u_2 - u_1)^n e^{(cst\ A_I\ \Delta t)}, \end{cases} \tag{A4}$$

with

$$cst = -\frac{Z_1 Z_2}{Z_1 + Z_2}\left(\frac{1}{(h_1\rho_1)^n} + \frac{1}{(h_2\rho_2)^n}\right). \tag{A5}$$

The conservative variables of the momentum equations are then updated after the hyperbolic step, the cylindrical correction step, and pressure relaxation as:

$$\begin{cases} (h_1\rho_1 u_1)^{n+1} = (h_1\rho_1)^{n+1} (u_1)_t, \\ (h_2\rho_2 u_2)^{n+1} = (h_2\rho_2)^{n+1} (u_2)_t, \end{cases} \tag{A6}$$

where $(u_1)_t$ and $(u_2)_t$ are given by (A4).

In the present contribution, the interfacial area has been determined in Section 5 as a function of time:

$$A_I(t) = a.h_0 e^{-b.t} + c, \tag{A7}$$

where $h_0$ is the initial height of the krypton column and $a, b, c$ are three positive constants.

However, System (A3) still admits explicitly solution as:

$$\begin{cases} (u_2)_t = \dfrac{(h_1\rho_1 u_1 + h_2\rho_2 u_2)^n + (h_1\rho_1)^n (u_2 - u_1)^n e^{cst\left(\frac{a.h_0}{-b} e^{-b.t_0}\left(e^{-b\Delta t} - 1\right) + c\Delta t\right)}}{(h_1\rho_1)^n + (h_2\rho_2)^n}, \\ (u_1)_t = (u_2)_t - (u_2 - u_1)^n e^{cst\left(\frac{a.h_0}{-b} e^{-b.t_0}\left(e^{-b\Delta t} - 1\right) + c\Delta t\right)}, \end{cases} \tag{A8}$$

with $t_0$ the current time. The conservative variables of the momentum equations are then updated with the help of Equation (A6).

## References

1. Brighton, P.W.M.; Prince, A.J.; Webber, D.M. Determination of cloud area and path from visual and concentration records. *J. Hazard. Mater.* **1985**, *11*, 155–178.
2. Deaves, D.M. Dense gas dispersion modelling. *J. Loss Prev. Process Ind.* **1992**, *5*, 219–227.
3. Ermak, D.L. *User's Manual for SLAB: An Atmospheric Dispersion Model for Denser-Than-Air Releases*; UCRL-MA-l05607, Lawrence Livermore National Laboratory, Livermore, C A. **1990**.
4. Hankin, R.K.S. Heavy Gas Dispersion over Complex Terrain. Ph.D. Thesis, University of Cambridge, UK, **1997**.
5. Hank, S.; Saurel, R.; LeMétayer, O.; Lapébie, E. Modeling blast waves, gas and particles dispersion in urban and hilly ground areas. *J. Hazard. Mater.* **2014**, *280*, 436–449.
6. Hank, S. Modélisation et Simulation de la Dispersion de Fluide en Milieu Fortement Hétérogène. Ph.D. Thesis, Aix Marseille Université, Marseille, France, **2012**.
7. Hanna, S.; Dharmavaram, S.; Zhang, J.; Sykes, I.; Witlox, H.; Khajehnajafi, S.; Koslan, K. Comparison of six widely-used dense gas dispersion models for three recent chlorine railcar accidents. *Process Saf. Prog.* **2008,** *27*, 248–259.
8. Chiapolino, A.; Saurel, R. Models and methods for two-layer shallow water flows. *J. Comput. Phys.* **2018**, *371*, 1043–1066.
9. Teshukov, V.M. Gas-Dynamic analogy for vortex free-boundary flows. *J. Appl. Mech. Tech. Phys.* **2007**, *48*, 303–309.
10. Richard, G.; Gavrilyuk, S. A new model of roll waves: Comparison with Brock's experiments. *J. Fluid Mech.* **2012**, *698*, 374–405.
11. Richard, G. Elaboration D'un Modèle D'écoulements Turbulents en Faible Profondeur. Application au Ressaut Hydraulique et aux Trains de Rouleaux. Ph.D. Thesis, Aix-Marseille Université, Marseille, France, **2013**.
12. Richard, G.; Gavrilyuk, S. The classical hydraulic jump in a model of shear shallow-water flows. *J. Fluid Mech.* **2013**, *725*, 492–521.
13. Richard, G.; Gavrilyuk, S. Modelling turbulence generation in solitary waves on shear shallow water flows. *J. Fluid Mech.* **2015**, *773*, 49–74.
14. Gavrilyuk, S.; Liapidevskii, V.; Chesnokov, A. Spilling breakers in shallow water: Applications to Favre waves and to the shoaling and breaking of solitary waves. *J. Fluid Mech.* **2016**, *808*, 441–468.
15. Saurel, R.; Pantano, C. Diffuse interface capturing methods for compressible two-phase flows. *Annu. Rev. Fluid Mech.* **2018**, *50*, 105–130.
16. Marble, F. Dynamics of a gas containing small solid particles. *Combustion and Propulsion (5th AGARD Colloquium), 175, Pergamon Press,* **1963**.
17. Baer, M.; Nunziato, J. A two-phase mixture theory for the deflagration-to-detonation transition (DDT) in reactive granular materials. *Int. J. Multiph. Flow* **1986**, *12*, 861–889.
18. Saurel, R.; Chinnayya, A.; Carmouze, Q. Modelling compressible dense and dilute two-phase flows. *Phys. Fluids* **2017**, *29*, 063301.
19. Lallemand, M.; Saurel, R. Pressure relaxation procedures for multiphase compressible flows. *Tech. Rep. INRIA Rep.* **2005**, doi:10.1002/fld.967.
20. Forestier, A.; Hérard, J.M.; Louis, X. Solveur de type Godunov pour simuler les écoulements turbulents compressibles. *Comptes Rendus l'Académie Sci. Paris Série Mathématique* **1997**, *324*, 919–926.
21. Saurel, R.; Gavrilyuk, S.; Renaud, F. A multiphase model with internal degrees of freedom: Application to shock–bubble interaction. *J. Fluid Mech.* **2003**, *495*, 283–321.
22. Lhuillier, D.; Chang, C., Theofanous, T. On the quest for a hyperbolic effective-field model of disperse flows. *J. Fluid Mech.* **2013**, *731*, 184–194.
23. Abgrall, R.; Karni, S. Two-Layer shallow water system: A relaxation approach. *SIAM J. Sci. Comput.* **2009**, *31*, 1603–1627.
24. Gallouet, T.; Masella, J. Un schéma de Godunov approché. *Comptes Rendus l'Académie Sci. Paris Série Mathématique* **1996**, *323*, 77–84.

25. Saurel, R.; Petitpas, F.; Berry, R. Simple and efficient relaxation methods for interfaces separating compressible fluids, cavitating flows and shocks in multiphase mixtures. *J. Comput. Phys.* **2009**, *228*, 1678–1712.
26. Kapila, A.; Menikoff, R., Bdzil, J.; Son, S.; Stewart, D. Two-Phase modeling of deflagration-to-detonation transition in granular materials: Reduced equations. *Phys. Fluids* **2001**, *13*, 3002–3024.
27. Chiapolino, A.; Saurel, R.; Nkonga, B. Sharpening diffuse interfaces with compressible fluids on unstructured meshes. *J. Comput. Phys.* **2017**, *340*, 389–417.
28. Ovsyannikov, L. Two-layer "shallow water" model. *J. Appl. Mech. Tech. Phys.* **1979**, *20*, 127–135.
29. Kurganov, A.; Petrova, G. Central-Upwind schemes for two-layer shallow water equations. *SIAM J. Sci. Comput.* **2009**, *31*, 1742–1773.
30. Monjarret, R. Local well-posedness of the two-layer shallow water model with free surface. *SIAM J. Appl. Math.* **2015**, *75*, 2311–2332.
31. Chinnayya, A.; Daniel, E.; Saurel, R. Modelling detonation waves in heterogeneous energetic materials. *J. Comput. Phys.* **2004**, *196*, 490–538.
32. Harten, A.; Lax, P.; van Leer, B. On upstream differencing and godunov-type schemes for hyperbolic conservation laws. *SIAM Rev.* **1983**, *25*, 35–61.
33. Davis, S. Simplified second-order Godunov-type methods. *SIAM J. Sci. Stat. Comput.* **1988**, *9*, 445–473.
34. Adduce, C.; Sciortino, G.; Proietti, S. Gravity currents produced by lock exchanges: Experiments and simulations with a two-layer shallow-water model with entrainment. *J. Hydraul. Eng.* **2012**, *198*, 111–121.
35. McQuaid, J.; Roebuck, B. Large scale field trials in dense vapour dispersion, final report to sponsors on the heavy gas dispersion trials at thorney island 1982–84. *Tech. Rep. Health Saf. Exec.* **1985**, doi:10.1007/978-94-009-4972-0_18.
36. Edgerton, H.E. Shock wave photography of large subjects in daylight. *Rev. Sci. Instrum.* **1958**, *29*, 171–172.
37. Van der Walt, S.; Schönberger, J. L.; Nunez-Iglesias, J.; Boulogne, F.; Warner, J.D.; Yager, N.; Gouillart, E.; Yu, T. Scikit-Image: Image processing in python. *Peer J.* **2014**, *2*, e453.
38. Frangi, A.F.; Niessen, W.J.; Vincken, K.L.; Viergever, M.A. Multiscale vessel enhancement filtering. *Int. Conf. Med Image Comput. Comput. Assist. Interv.* **1998**, *1496*, 130–137.
39. Canny, J. A computational approach to edge detection. *IEEE Trans. Pattern Anal. Mach. Intell.* **1986**, *8*, 679–698.
40. Nielsen, M.; Søren, O. *A Collection of Data from Dense Gas Experiments*. Technical report Risø-R-845, Risø National Laboratory: Risoe-R, Forskningscenter Risoe, Denmark, **1995**.
41. Lord Kelvin Hydrokinetic solutions and observations. *Philos. Mag.* **1871**, *42*, 362–377.
42. Helmholtz, H. Űber discontinuirliche Flüssigkeitsbewegungen. *Wiss. Abh.* **1882**, 146–157.
43. Lord Rayleigh Investigation of the character of the equilibrium of an incompressible heavy fluid of variable density. *Sci. Pap.* **1900**, *2*, 200–207.
44. Taylor, G.I. The instability of liquid surfaces when accelerated in a direction perpendicular to their planes. I. *Proc. R. Soc. Lond. Ser.* **1950**, *201*, 192–196.
45. van Leer, B. Towards the ultimate conservative difference scheme. II. Monotonicity and conservation combined in a second-order scheme. *J. Comput. Phys.* **1974**, *14*, 361–370.